
\input amstex
\documentstyle{amsppt}
\nologo
\catcode`\@=11
\footline={\hfil}
\font\bigtitlefont=cmr17
\font\bigauthorfont=cmr12
\def\title{\let\savedef@\title
  \def\title##1\endtitle{\let\title\savedef@
    \global\setbox\titlebox@\vtop{\bigtitlefont
      \raggedcenter@
      \baselineskip1.3\baselineskip{##1}\endgraf}%
    \ifmonograph@
      \edef\next{\the\leftheadtoks}\ifx\next\empty@ \leftheadtext{##1}\fi
    \fi
    \edef\next{\the\rightheadtoks}\ifx\next\empty@ \rightheadtext{##1}\fi
  }%
  \nofrillscheck\title}
\def\author#1\endauthor{\global\setbox\authorbox@
  \vbox{\bigauthorfont\raggedcenter@ #1\endgraf}\relaxnext@
  \edef\next{\the\leftheadtoks}%
  \ifx\next\empty@\leftheadtext{#1}\fi}
\def\address#1\endaddress{\global\advance\addresscount@\@ne
  \expandafter\gdef\csname address\number\addresscount@\endcsname
  {\nobreak\vskip12\p@ minus6\p@\noindent\eightpoint\smc#1\par}}
\def\email{\let\savedef@\email
  \def\email##1\endemail{\let\email\savedef@
  \toks@{\def\usualspace{{\it\enspace}}\endgraf\noindent\eightpoint}%
  \toks@@{\tt##1\par}%
  \expandafter\xdef\csname email\number\addresscount@\endcsname
  {\the\toks@\frills@{{\noexpand\it E-mail address\noexpand\/}:%
     \noexpand\enspace}\the\toks@@}}%
  \nofrillscheck\email}
\font@\refsheadfont@=cmbx10 at 14pt
\def\Refsname{References\kern3.4in}
\def\refsfont@{\tenpoint}
\catcode`\@=\active

\pageno=165
\hoffset=1in
\voffset=.75in
\pagewidth{27pc}
\pageheight{520bp}
\baselineskip=13pt

\def\folio{\number\pageno}
\headline={\ifnum\pageno=165{\hfil}\else
  {\ifodd\pageno{{\tensl Automorphisms of manifolds}\hfil{\tenrm\folio}}
   \else{{\tenrm\folio}\hfil{\tensl Michael Weiss and Bruce Williams}}
\fi}\fi}

\loadbold
\loadeusm
\define\la{\longrightarrow}

\define\minus{\smallsetminus}

\define\ep{\varepsilon}

\define\bo{\bold b\bold o}
\define\bhm{\bold b\bold h\bold m}
\define\td{^{\sim}}
\define\wdt{\widetilde}

\redefine\tan{_{\text{tan}}}
\define\itb{\item"{$\bullet$}"}
\define\upr{^{\%}}
\define\lpr{_{\%}}

\define\sha{^{\sharp}}

\define\nat{^{\natural}}
\define\8{^{\infty}}
\define\hf{^{h\ZZ/2}}
\define\ho{_{h\ZZ/2}}
\define\thf{^{th\ZZ/2}}

\define\bul{_{\bullet}}
\define\ubul{^{\bullet}}

\define\sB{\eusm B}
\define\sC{\eusm C}
\define\sD{\eusm D}

\define\sH{\eusm H}

\define\sR{\eusm R}
\define\sS{\eusm S}
\define\sT{\eusm T}

\define\sW{\eusm W}

\define\BB{\Bbb B}
\define\CC{\Bbb C}
\define\DD{\Bbb D}

\define\HH{\Bbb H}

\define\NN{\Bbb N}
\define\QQ{\Bbb Q}
\define\RR{\Bbb R}
\define\SS{\Bbb S}
\define\ZZ{\Bbb Z}

\define\bA{\bold A}
\define\bC{\bold C}

\define\bF{\bold F}

\define\bH{\bold H}
\define\bJ{\bold J}
\define\bK{\bold K}
\define\bL{\bold L}
\define\bS{\bold S}
\define\bT{\bold T}

\define\bV{\bold V}

\define\bY{\bold Y}

\define\id{\operatorname{id}}

\define\hocolim{\operatornamewithlimits{hocolim}}
\define\hofiber{\operatorname{hofiber}}

\define\im{\operatorname{im}}
\define\Whd{\bold W\bold h\bold d}

\define\TOP{\operatorname{TOP}}
\define\DIFF{\operatorname{DIFF}}
\define\G{\operatorname{G}}
\define\C{\operatorname{C}}
\define\map{\operatorname{map}}
\define\Or{\operatorname{O}}

\define\Nil{\text{{\bf Nil}}}
\define\msk{\medskip\noindent}
\define\bsk{\bigskip\noindent}

\topmatter
\title Automorphisms of manifolds 
\endtitle
\author Michael Weiss and Bruce Williams 
\endauthor
\address Dept. of Mathematics\newline
Univ. of Aberdeen\newline
Aberdeen AB24 3UE, U\.K\.
\endaddress
\email m.weiss\@maths.abdn.ac.uk \endemail
\address Dept\. of Mathematics\newline
University of Notre Dame\newline
 Notre Dame, IN 46556, USA
\endaddress 
\email williams.4\@nd.edu \endemail 
\thanks Both authors partially supported by NSF grants. \endthanks
\endtopmatter

\head {\bf 0. Introduction} \endhead 
This survey is about homotopy types of spaces of automorphisms of
topological  and smooth manifolds. Most of the results available are 
{\it relative}, i\.e\., they compare different types of automorphisms.

In chapter 1, which motivates the later chapters, we introduce our
favorite types of  manifold automorphisms and 
make a comparison by (mostly elementary) geometric methods. 
Chapters 2, 3, and 4 describe algebraic models (involving $L$--theory
and/or algebraic $K$--theory) for certain spaces of ``structures''
associated  with a manifold $M$, that is, spaces of other manifolds 
sharing some geometric features with $M$. 
The algebraic models rely heavily on
\roster
\itb Wall's work in surgery theory, e\.g\.  \cite{Wa1}~,
\itb Waldhausen's
work in $h$--cobordism  theory alias concordance theory, 
which includes a parametrized  version of Wall's theory of the finiteness
obstruction, \cite{Wa2}\,. 
\endroster
The structure spaces are of interest for the  following reason. Suppose that
two different notions of automorphism of $M$ are being compared.
Let $X_1(M)$ and $X_2(M)$ be the
corresponding automorphism  spaces; suppose that 
$X_1(M)\subset X_2(M)$. As a  rule, the space of cosets
$X_2(M)/X_1(M)$ is a union of connected  components of a suitable
structure space. 

Chapter 5 contains the
beginnings of a more radical approach in which one tries to  calculate
the classifying space
$BX_1(M)$ in terms of $BX_2(M)$,  rather  than trying to calculate
$X_2(M)/X_1(M)$.  Chapter 6 contains some
examples and calculations. 

Not included in this survey is the  
disjunction theory of automorphism spaces and embedding spaces begun by 
Morlet, see \cite{BLR}, and continued in 
\cite{Go4},  \cite{Go1}, \cite{Go2}, \cite{We2}, \cite{We3}, \cite{GoWe},
\cite{GoKl}. It calls for a survey of its own.
\bsk
\head {\bf 1. Stabilization and descent} \endhead
\head 1.1.  Notation, terminology \endhead
\definition{1.1.1. Terminology} {\it Space} with a capital $S$ 
means {\it simplicial set}. We will occasionally see simplicial Spaces
(=bisimplicial sets). A  simplicial Space $k\mapsto Z_k$ determines a Space 
$(\amalg_k\,\,\Delta^k\times Z_k)/\sim$~,
where $\Delta^k$ is the $k$--simplex viewed as a Space (= simplicial
set) and $\sim$ refers to the relations $(f_*x,y)\sim(x,f^*y)$. The Space
$(\amalg_k\,\,\Delta^k\times Z_k)/\sim$ is isomorphic to the Space
$k\mapsto Z_k(k)$, the {\it diagonal} of
$k\mapsto Z_k$.  See  \cite{Qui}, for example.

A {\it euclidean $k$--bundle} is a fiber bundle with 
fibers homeomorphic to $\RR^k$. Trivial euclidean $k$--bundles 
are often denoted $\ep^k$.  

The homotopy fiber of a map $B\to C$ of Spaces, where $C$ is 
based, will be denoted $\hofiber[\,B\to C\,]$. A {\it homotopy 
fiber sequence} is a diagram of spaces 
$A@>>>B@>>>C$ where $C$ is based, together with a 
nullhomotopy of the composition $A\to C$ which makes the 
resulting map from $A$ to $\hofiber[\,B\to C\,]$ a (weak) homotopy 
equivalence. 

The term {\it cartesian square} is synonymous with {\it homotopy pullback 
square}. More generally, an {\it $n$--cartesian square} ($n\le\infty)$ 
is a commutative diagram of Spaces and maps
$$\CD
A@>>> C\\
@VVV @VVV \\
B @>>> D 
\endCD$$
such that the resulting map from $A$ to the homotopy
pullback of the diagram $B@>>>D@<<<C$ is $n$--connected.  

A {\it commutative diagram} of Spaces is a functor $F$ 
from some small category $\sD$ to Spaces. We say that $\sD$ 
is the shape of the diagram. When we
represent  such a  diagram graphically, 
we usually only show the maps 
$F(g_i)$ for a set $\{g_i\}$ of morphisms generating $\sD$. For 
example, the commutative square just above is a functor from 
a category with four objects and five non--identity morphisms to
Spaces.  The notion of a {\it homotopy commutative diagram} of
shape $\sD$ has been made precise by \cite{Vogt}. It is a 
continuous functor from a certain topological category $\sW\sD$
(determined  by $\sD$) to Spaces. In more detail, 
$\sW\sD$ is a small topological category with discrete object set, and
comes with a continuous functor $\sW\sD\to \sD$ which restricts to
a  homeomorphism (=bijection) of object sets, and to a homotopy 
equivalence of morphism spaces. Graphically, we represent 
homotopy commutative diagrams of shape $\sD$ like 
commutative diagrams of shape $\sD$. 
\enddefinition
\msk
\definition{1.1.2.  Notation} For a topological 
manifold $M$, we denote by $\TOP(M)$ the Space of homeomorphisms
$f:M\to M$ which agree with the identity on $\partial M$. (A
$k$--simplex in $\TOP(M)$ is a homeomorphism $f:M\times\Delta^k\to
M\times\Delta^k$   over $\Delta^k$ which agrees with the identity on $\partial
M\times\Delta^k$.) References: \cite{BLR}, \cite{Bu1}.

We use the abbreviations $\TOP(n)=\TOP(\RR^n)$ and 
$\TOP=\bigcup_n\TOP(n)$, where we include $\TOP(n)$ in 
$\TOP(n+1)$ by $f\mapsto f\times\id_{\RR}$. 

Let $\partial_+M$ be a codimension zero 
submanifold of $\partial M$, closed as a subspace of $\partial M$.
Let $\partial_-M$ be the closure of $\partial M\minus\partial_+M$. 
Let  $\TOP(M,\partial_+M)$ be 
the Space of homeomorphisms $M\to M$ which agree
with the identity on $\partial_-M$.  (This is $\TOP(M)$ 
if $\partial_+M=\emptyset$.) Special case:  the Space of 
concordances alias pseudo--isotopies of a 
compact manifold $N$, which is 
$\C(N):=\TOP(N\times I,N\times 1)$ where $I=[0,1]$.  References:
\cite{Ce}, \cite{HaWa}, \cite{Ha}, \cite{DIg}, \cite{Ig}. 
\comment
When $M$ is smooth and $\partial_+M$ is a smooth submanifold of 
$\partial M$, then $\DIFF(M,\partial_+M)$ is defined; it is the  
Space of diffeomorphisms $M\to M$ (which agree with the identity on 
$\partial_-M$.  Also, there is a simplicial group 
$\C_d(M)$ of {\it smooth} concordances.
\endcomment
Let $\G(M,\partial_+ M)$ 
be the Space of homotopy equivalences of triads,
$(M;\partial_+M,\partial_-M)\to (M;\partial_+M,\partial_-M)$, 
which are the identity on $\partial_-M$.  We 
abbreviate $\G(M,\emptyset)$ to $\G(M)$.  {\it Warning:} If
$\partial M=\emptyset$, then $\G(M)$  is the Space of 
homotopy equivalences  $M\to M$, but in general it is not. 
\enddefinition
\msk
\definition{1.1.3. Definitions} An {\it $h$--structure} on a closed manifold
$M^n$ is a pair $(N,f)$ where $N^n$ is another closed manifold and
$f:N\to M$ is a homotopy equivalence.  If $f$ is a 
{\it simple} homotopy equivalence, $(N,f)$ is an
{\it $s$--structure}. An {\it isomorphism} 
from an $h$--structure $(N_1,f_1)$ to another $h$--structure $(N_2,f_2)$ on 
$M$ is a homeomorphism $N_1\to N_2$ over $M$. 

We see that the $h$--structures on $M$ form a groupoid. Better,
they form a {\it simplicial groupoid}:  Objects in degree $k$ are 
pairs $(N,f)$ where $N^n$ is another closed manifold 
and $f:N\times\Delta^k\to M\times\Delta^k$ is a homotopy 
equivalence over $\Delta^k$. Morphisms in degree $k$ are 
homeomorphisms over $M\times \Delta^k$.

Let $\sS(M)$ be the diagonal nerve (= diagonal of degreewise nerve)
of this simplicial groupoid; also let $\sS^s(M)$ be the diagonal 
nerve of the simplicial subgroupoid of $s$--structures.
Think of $\sS(M)$ and $\sS^s(M)$ as the {\it Spaces} of 
$h$--structures on $M$ and $s$--structures on $M$, respectively.  (They
are actually simplicial {\it classes}, not  simplicial sets, as it stands. The
reader can either accept this, or avoid it by working in a Grothendieck
``universe''.)  The forgetful  functor $(N,f)\mapsto N$ induces  a map from
$\sS(M)$ to  the diagonal nerve  of  the simplicial groupoid of 
all closed $n$--manifolds and homeomorphisms between such.
This map is a Kan fibration. 
Its fiber over the point corresponding to $M$ is $\G(M)$.  Hence 
there is a homotopy fiber sequence 
$$\TOP(M)@>>>\G(M)@>>>\sS(M)\,.$$
More generally, given compact $M$ and $\partial_+M\subset\partial M$
as above, there is an $h$--structure Space $\sS(M,\partial_+M)$ and 
an $s$--structure Space $\sS^s(M,\partial_+M)$ 
and a homotopy fiber sequence 
$$\TOP(M,\partial_+M)@>>>\G(M,\partial_+M)@>>>\sS(M,\partial_+M)\,.$$
We omit the details. Important special cases: the 
{\it Space of $h$--cobordisms} and the {\it Space of $s$--cobordisms}
on a compact manifold,
$$\align \sH(N)&:=\sS(N\times I,N\times 1)\\
\sH^s(N)&:=\sS^s(N\times I,N\times 1)\,.\endalign$$
 Since $\G(N\times I,N\times 1)\,\simeq\,*$, we have 
$\Omega\sH(N)\simeq \C(N)$. 
There is a stabilization map 
$\sH(N)\to\sH(N\times I)$, {\it upper stabilization} to be precise
\cite{Wah2}, \cite{HaWa}.  Let 
$$\sH^{\infty}(N):=\hocolim_k \sH(N\times I^k)\,,\qquad
\C^{\infty}(N):=\hocolim_k \C(N\times I^k)\,.$$
\enddefinition
\msk
\definition{1.1.4. More definitions} 
The {\it block automorphism Space}
$\,\wdt{\TOP}(M)$ has as its $k$--simplices 
the homeomorphisms 
$g:M\times\Delta^k\to M\times\Delta^k$ which satisfy
$g(M\times s)=M\times s$ for each face $s\subset\Delta^k$, 
and restrict to the identity on $\partial M\times\Delta^k$. References: \cite{ABK}, \cite{Bu1}
\cite{Br1}. There is also a block $s$-structure Space 
$$\wdt{\sS}^s(M),$$ 
defined as the diagonal nerve of a simplicial groupoid. The objects 
of the simplicial groupoid in degree $k$ are of the form 
$(N,f)$ where $N^n$ is closed and $f$ is a {\it simple} homotopy equivalence
$N\times\Delta^k\to M\times\Delta^k$ 
such that $f(N\times t)\subset M\times t$ 
for each face $t$ of $\Delta^k$, and $f$ restricts to a
homeomorphism 
$\partial N\times\Delta^k\to \partial M\times\Delta^k$. The morphisms 
in degree $k$ are homeomorphisms respecting the reference maps to
$M\times\Delta^k$. References: \cite{Qun1}, 
\cite{Wa1, \S17.A}, \cite{Ni}, \cite{Rou1}.
There is a homotopy fiber sequence 
$$\wdt{\TOP}(M)@>>>\wdt{\G}^s(M)@>>>\wdt{\sS}^s(M),$$
where $\wdt{\G}^s(M)$ is defined like $\wdt{\TOP}(M)$, but with 
simple homotopy equivalences instead of homeomorphisms.--- 
These definitions have relative versions (details omitted); for 
example, there is a homotopy fiber sequence 

$$\wdt{\TOP}(M,\partial_+M)@>>>\wdt{\G}^s(M,\partial_+M)
@>>>\wdt{\sS}^s(M,\partial_+M)\,.$$

Let $\G^s(M,\partial_+M)\subset \G(M,\partial_+M)$ consist 
of the components containing those $f$ which are simple 
homotopy automorphisms and induce simple homotopy 
automorphisms of $\partial_+M$. 
The inclusion 
$$\G^s(M,\partial_+M)\to \wdt{\G}^s(M,\partial_+M)$$
 is a homotopy equivalence (because it induces an isomorphism on 
homotopy groups;  both Spaces are fibrant).
\enddefinition
\bsk
\head 1.2.  Open stabilization versus closed stabilization \endhead
Let $M^n$ be compact, $M_0=M\minus\partial M$.
{\it Open stabilization} refers to the map
$$\TOP(M,\partial M)@>>>\bigcup_k\TOP(M_0\times\RR^k)$$
given by $f\mapsto f|M_0$. We include $\TOP(M_0\times\RR^k)$ 
in $\TOP(M_0\times\RR^{k+1})$ by $g\mapsto g\times\id_{\RR}$. 
{\it Closed stabilization} refers to the inclusion 
$$\TOP(M,\partial M)@>>>\bigcup_k\TOP(M\times I^k,\partial(M\times I^k))\,.$$
Open stabilization factors through closed stabilization, by 
means of the restriction maps  
$\TOP(M\times I^k,\partial(M\times I^k))\to 
\TOP(M_0\times I^k_0)$ and an identification $I_0\cong\RR$. 
Here in \S1.2 we describe the homotopy type of 
$\bigcup_k\TOP(M_0\times\RR^k)$, and descend 
from there to  $\bigcup_k\TOP(M\times I^k,\partial(M\times I^k))$. 
For a more algebraic version of this, see \S5.2. 
\bsk
Let $\hat\tau:M_0\to B\!\TOP(n)$ classify the tangent bundle \cite{Mi1}, \cite{Kis}, 
\cite{Maz2}, \cite{KiSi,IV.1}. We  
map $\G(M_0)$ to the mapping Space $\map(M_0, B\!\TOP)$ by 
$f\mapsto\hat\tau f$
and we map $\bigcup_k\TOP(M_0\times\RR^k)$ to  $\G(M_0)$ 
by $f\mapsto pfi$, where $p:M_0\times\RR^k\to  M_0$ and 
$i:M_0\cong M_0\times \bold 0\to M_0\times\RR^k$ are projection 
and inclusion, respectively. 
\msk
\proclaim{1.2.1. Theorem \cite{CaGo}} The resulting diagram 
$$\bigcup_k\TOP(M_0\times\RR^k)\la \G(M_0) \la \map(M_0, B\!\TOP)$$
is a homotopy fiber sequence.
\endproclaim
\msk
The proof uses immersion theory \cite{Gau},  
general position, and the half--open $s$--cobordism theorem 
\cite{Sta}.  See also \cite{Maz1}. 
\msk
Choose a collar for $M$, that is, an embedding $c:\partial M\times I\to M$
extending the map $(x,0)\to x$ on $\partial M\times 0$. Reference: \cite{Brn}, 
\cite{KiSi, I App\. A}. Any 
homeomorphism $f:M_0\to M_0$ determines an $h$--cobordism $W_f$
on $\partial M$: the region of $M$ enclosed by $\partial M$ 
and $fc(\partial M\times 1)$. The bundle on the geometric 
realization of $\TOP(M_0)$ with fiber $W_f$ over the vertex $f$ is a
bundle  of $h$--cobordisms, classified by a map $v$ from
$\TOP(M_0)$ to $\sH(\partial M)$.

If $f:M_0\to M_0$ is the restriction of 
some homeomorphism $g:M\to M$, then 
$W_f\cong gc(\partial M\times I)$ is trivialized. 
Conversely, a trivialization of $W_f$ can be used 
to construct a homeomorphism $g:M\to M$ with an isotopy from 
$g|M_0$ to $f$. Therefore: the diagram
$$\TOP(M,\partial M)@>\text{res}>>\TOP(M_0)@>v>>\sH(\partial M)
\tag1.2.2$$
is a homotopy fiber sequence. See \cite{Cm} for details. The special case 
where $M=\DD^n$ is due to \cite{KuLa}. 

This observation can be stabilized. Let $u:\sH(\partial M)\to 
\sH(\partial(M\times I))$ be the composition of stabilization 
$\sH(\partial M)\to \sH(\partial M\times I)$ with the map 
induced by the inclusion of $\partial M\times I$ in $\partial(M\times I)$.
Then 
$$\CD
\TOP(M_0)@>v>>\sH(\partial M) \\
@VVV   @VV u V  \\
\TOP(M_0\times I_0)@>>>\sH(\partial(M\times I))
\endCD$$
is homotopy commutative.   The homotopy colimit of the 
$\sH(\partial(M\times I^k))$ under the $u $--maps 
becomes $\simeq \sH^{\infty}(M)$. Therefore: 
\msk
\proclaim{1.2.3. Theorem} There exists a homotopy fiber sequence
$$\bigcup_k\TOP(M\times I^k,\partial(M\times I^k))
@>{\quad\text{res}\quad}>>\bigcup_k\TOP(M_0\times\RR^k) 
@>\qquad>>\sH^{\infty}(M)\,.$$
\endproclaim
\msk
Let $Q=I\8$ be the Hilbert cube. The product $M\times Q$ 
is a {\it Hilbert cube manifold} \cite{Cha}  without boundary.
Let $\TOP(M\times Q)$ be the Space of homeomorphisms 
$M\times Q\to M\times Q$ and let $\G(M\times Q)$ be the 
Space of homotopy equivalences $M\times Q\to M\times Q$. 
Chapman and Ferry have shown \cite{Bu2} that an evident map
from $\bigcup_k\TOP(M\times I^k,\partial(M\times I^k))$ to the 
homotopy fiber of the composition 
$$\TOP(M\times Q)@>\subset>>\G(M\times Q)\,\,\simeq\,\,
G(M_0)@>>>\map(M_0,B\!\TOP)$$
(last arrow as in 1.2.1) is a homotopy 
equivalence. Therefore 
$$\bigcup_k\TOP(M\times I^k,\partial(M\times I^k))
\la \TOP(M\times Q)\la \map(M_0,B\!\TOP)$$
is a homotopy fiber sequence.  Comparison  with 1.2.1
gives  the next result. 
\msk 
\proclaim{1.2.4. Theorem} The following homotopy 
commutative diagram is
cartesian:
$$\CD
\bigcup_k\TOP((M\times I^k,\partial(M\times I^k))
@>\text{res}>>\bigcup_k\TOP(M_0\times\RR^k) \\
@VVV      @VVV \\
\TOP(M\times Q) @>>> \G(M_0)\,.
\endCD$$
\endproclaim
\msk
This suggests that the map 
$\bigcup_k\TOP(M_0\times\RR^k)\to\sH^{\infty}(M)$ in 1.2.3 
factors through $\G(M_0)$.  We will obtain such a factorization in 1.5.3. 
\msk
\remark{Remark} 
Looking at horizontal homotopy fibers in 1.2.4, and using 
1.2.3, and the homotopy equivalence $\G(M_0)\simeq \G(M\times Q)$,
one finds that the homotopy fiber of 
the inclusion $\TOP(M\times Q)\to \G(M\times Q)$ is 
$\C^{\infty}(M)$. This can also be deduced from  \cite{Cha2}, \cite{Cha3}. 
\endremark
\bsk
\head  1.3.  Bounded stabilization versus no stabilization \endhead
Let $M^n$ be compact. A homeomorphism 
$f:M\times\RR^k\to M\times\RR^k$ is {\it bounded} if $\{p_2f(z)-p_2(z)\mid z\in M\times\RR^k\}$ is a 
bounded subset of $\RR^k$, where $p_2:M\times\RR^k\to\RR^k$ is the projection.
Let $\TOP^b(M\times\RR^k)$ be the Space of bounded
homeomorphisms $M\times\RR^k\to M\times\RR^k$ which agree
with the identity on $\partial M\times\RR^k$. 
Note $\TOP(M)=\TOP^b(M\times\RR^0)$. 
{\it Bounded stabilization} refers to the inclusion 
$$\TOP(M)@>>> \bigcup_k \TOP^b(M\times\RR^k)\,.$$
Surgery theory describes the homotopy type of $\bigcup_k
\TOP^b(M\times\RR^k)$, modulo the mysteries of $\G(M)$. See \S2.4; here
in \S1.3 we analyze the difference between  $\bigcup_k
\TOP^b(M\times\RR^k)$ and 
$\TOP(M)$. 
\bsk
The Space $\TOP^b(M\times\RR^{k+1})/\TOP^b(M\times\RR^k)$
for $k\ge 0$ and fixed $M$ is the $k$--th Space in a spectrum $\bH(M)$, by analogy with the 
sphere spectrum, which is made out of the spaces 
$\Or(\RR^{k+1})/\Or(\RR^k)$. Compare \cite{BuLa1}. Anderson and Hsiang, who 
introduced bounded homeomorphisms in \cite{AH1}, \cite{AH2} showed that 
$\Omega^{\infty+1}(\bH(M))\simeq \C^{\infty}(M)$.  In more detail: they
introduced {\it bounded concordance spaces} 
$$\C^b(M\times\RR^k)=
\TOP^b(M\times I\times \RR^k,M\times 1\times\RR^k)$$
and proved the following. See also \cite{WW1,\S1+App.5}, \cite{Ha, App.II}.  
\msk
\proclaim{1.3.1. Theorem \cite{AH1}, \cite{AH2}} Assume $n>4$. Then 
\roster
\item"{i)}" $\Omega(\TOP^b(M\times\RR^{k+1})/\TOP^b(M\times\RR^k))
\,\simeq\, \C^b(M\times\RR^k)$;
\item"{ii)}" $\Omega\C^b(M\times\RR^k)\,\simeq\, 
\C^b(M\times I\times\RR^{k-1})$.
\endroster
\endproclaim
\msk
Part ii) of 1.3.1 shows that the spaces $\C^b(M\times\RR^k)$ for $k\ge0$
form a spectrum, with structure maps
$$\C^b(M\times\RR^{k-1})@>\text{stab\.}>>
\C^b(M\times I\times\RR^{k-1})
\,\,\simeq\,\,\Omega\C^b(M\times\RR^k);$$
then part i) of 1.3.1 with some extra work \cite{WW1,\S1} identifies 
the new spectrum with $\Omega\bH(M)$. It is also shown in \cite{AH1} that 
the homotopy groups $\pi_j\bH(M)$ for $j\le0$ are lower $K$--groups
\cite{Ba}:
\msk
\proclaim{1.3.2. Theorem}  Let $j \leq 0$ be an integer.  Then

$$ \pi_j\C^b(M\times\RR^k)=\left\{\aligned
    & K_{j-k+2}(\ZZ\pi_1(M)) \qquad\qquad( j<k-2) \\
    & \wdt{K}_0(\ZZ\pi_1(M)) \qquad\qquad\qquad ( j=k-2) \\
    & \text{Whitehead gp\. of } \pi_1(M)\quad( j=k-1)\,.
\endaligned\right.$$
\endproclaim
\msk
\remark{Remark} Madsen and Rothenberg \cite{MaRo1} , \cite{MaRo2} have proved 
equivariant analogs of 1.3.1 and 1.3.2, and Chapman \cite{Cha4}, 
Hughes \cite{Hu} 
have proved a Hilbert cube analog. Carter \cite{Ca1}, \cite{Ca2}, 
\cite{Ca3} has shown that 
$K_r(\ZZ\pi)$ vanishes if $\pi$ is finite and $r<-1$. 
\endremark
\msk
\remark{Remark} It is shown in \cite{WW1,\S5} that $\Omega\8\bH(M)\simeq\sH\8(M)$; this 
improves slightly on $\Omega^{\infty+1}\bH(M)\simeq\C\8(M)$.
\endremark
\msk
Theorems 1.3.1 and 1.3.2 are about
descent from $\TOP^b(M\times\RR^{k+1})$ to
$\TOP^b(M\times\RR^k)$.  For instant descent from
$\TOP^b(M\times\RR^{k+1})$ to $\TOP(M)$, there is the
{\it hyperplane test} \cite{WW1,\S3}, \cite{We4}. Think of $\RR P^k$ as 
the Grassmannian of codimension one linear subspaces 
$W\subset \RR^{k+1}$.  Let $\Gamma_k$ be the 
Space of sections of the bundle $E(k)\to\RR P^k$ with fibers 
$$E(k)_W:=\TOP^b(M\times\RR^{k+1})/\TOP^b(M\times W)$$
(see the remark just below). Note that $E(k)\to\RR P^k$ has a trivial section picking the coset 
$[\id]$ in each fiber; so $\Gamma_k$ is a based Space. 
\msk
\remark{Remark} The ``bundle'' $E(k)\to\RR P^{k+1}$ is really 
a {\it twisted cartesian product} \cite{Cu} with base Space 
equal to the singular simplicial set of $\RR P^k$, and with fibers 
$E(k)_W$ over a vertex $W$ as stated.
\endremark
\msk
We define a map
$\Phi_k:\TOP^b(M\times\RR^{k+1})/\TOP(M)\to \Gamma_k$
by taking the coset $f\cdot\TOP(M)$ to the section 
$W\mapsto f\cdot\TOP^b(M\times W)$.  For $k>0$, it is easy to 
produce an embedding $v_k$ making the square
$$\CD
\TOP^b(M\times\RR^k)/\TOP(M) @>\Phi_{k-1}>> \Gamma_{k-1} \\
@V\cap VV    @V v_k VV \\
\TOP^b(M\times\RR^{k+1})/\TOP(M) @>\Phi_k>> \Gamma_k 
\endCD$$
commutative.  Let $\Phi:\bigcup_k\TOP^b(M\times\RR^k)/\TOP(M)
\to\bigcup_k\Gamma_k$ be the union of the $\Phi_{k-1}$ for $k\ge0$. It
turns  out that $\Phi$ is highly connected (1.3.5 below), under mild 
conditions on $M$.
\msk
\definition{1.3.3. Definition} An integer 
$j$ is in the {\it topological, resp\. smooth, concordance stable range} for 
$M$ if the upper stabilization maps from 
$\C(M\times I^r)$ to $\C(M\times I^{r+1})$, resp\. the smooth versions, 
are $j$--connected, for all $r\ge0$.
\enddefinition
\msk
\proclaim{1.3.4. Theorem \cite{Ig}} If $M$ is smooth and 
$n\ge\max\{2j+7,3j+4\}$, 
then $j$ is in the smooth and in the topological 
concordance stable range for $M$.  {\rm (The  
estimate for the topological concordance stable range is due
to Burghelea--Lashof  and Goodwillie. Their argument 
uses smoothness of $M$, and Igusa's estimate of the 
smooth concordance stable range.  See \cite{Ig,
Intro\.}\,)} 
\endproclaim
\msk
\proclaim{1.3.5. Proposition \cite{WW1}} If $j$ is in the topological
concordance  stable range for $M$, and $n>4$, then 
$\Phi:\bigcup_k\TOP^b(M\times\RR^k)/\TOP(M)\to 
\bigcup_k\Gamma_{k-1}$ is $(j+1)$--connected. 
\endproclaim
\msk
\demo{Outline of proof } For $-1\le\ell\le k$ let
$\Gamma_{k,\ell}\subset\Gamma_k$ consist of the sections $s$ for 
which $s(W)=*$ whenever $W$ contains the standard copy of
$\RR^{\ell+1}$ in $\RR^{k+1}$.
Let $\Phi_{k,\ell}$ be the restriction of $\Phi_k$ 
to $\TOP^b(M\times\RR^{\ell+1})/\TOP(M)$, viewed as a map with 
codomain $\Gamma_{k,\ell}$.  One shows by induction on $\ell$
that the $\Phi_{k,\ell}$ for fixed $\ell$ define a highly connected map
$$\TOP^b(M\times\RR^{\ell+1})/\TOP(M) \la \bigcup_{k\ge\ell}
\Gamma_{k,\ell}\,.\qed $$
\enddemo
\msk
\proclaim{1.3.6 Theorem \cite{WW1}} 
 There exists a homotopy equivalence $$\bigcup_k\Gamma_k\,\simeq\,
\Omega\8(\bH(M)\ho)$$
for some involution on $\bH(M)$. Here
$\bH(M)\ho:=(E\ZZ/2)_+\wedge_{\ZZ/2}\bH(M)$ is the homotopy 
orbit spectrum. 
\endproclaim
\msk
\demo{Sketch proof} Note $\bigcup_k\Gamma_k=
\bigcup_{\ell>0}\bigcup_{k>\ell}\Gamma_{k,\ell}$ and 
$\bigcup_{k>\ell}\Gamma_{k,\ell}$ is homotopy equivalent to 
$\Omega\8(\SS^{\ell}_+\wedge_{\ZZ/2}\bH(M))$ by Poincar\'e duality
\cite{WW1, 2.4}, using 1.3.1 and 1.3.2. \qed
\enddemo
\msk
\definition{1.3.7. Summary} There exist a spectrum $\bH(M)$ 
with involution and a $(j+1)$--connected map 
$$\bigcup_k\,\TOP^b(M\times\RR^k)/\TOP(M)\,\,\la\,\,
\Omega\8(\bH(M)\ho)$$
where $j$ is the largest integer in the topological concordance
stable range for $M$. Further,
$\Omega\8\bH(M)\,\simeq\,\sH\8(M)$ and 
$$\pi_r\bH(M):=\left\{\aligned &\text{Wh}_1(\pi_1(M))\qquad\qquad r=0 \\
&\wdt{K}_0(\ZZ\pi_1(M))\qquad\quad r=-1 \\
&K_{r+1}(\ZZ\pi_1(M))\qquad r<-1\,. \endaligned\right. $$
(See \S3 for information about $\pi_r\bH(M)$ when $r>0$.)
\enddefinition
\bsk
\head 1.4. Block automorphism Spaces \endhead
The property $\pi_k\wdt{\TOP}(M)\,\cong\,
\pi_0\wdt{\TOP}(M\times\Delta^k)$ is a consequence of the definitions 
and has made the block automorphism Spaces popular. 
See also \S2.  In homotopy theory terms, 
the block  automorphism Space of $M$ is more closely related 
to $\bigcup_k\TOP^b(M\times\RR^k)$ than to $\TOP(M)$. To
explain this we need the bounded block automorphism Spaces 
$$\wdt{\TOP}^b(M\times\RR^k)$$
(definition left to the reader). The following {\it Rothenberg} type 
sequence is obtained by inspection, using 1.3.2.  For notation, 
see 1.3.7. Compare \cite{Sha}, \cite{Ra1, \S1.10}. 
\msk
\proclaim{1.4.2. Proposition \cite{AnPe}, \cite{WW1}} For $k\ge0$ there 
exists a long exact sequence 
$$\align \cdots&@>>>\pi_r\wdt{\TOP}^b(M\!\times\!\RR^k)
@>>>\pi_r\wdt{\TOP}^b(M\!\times\!\RR^{k+1})
@>>>H_{r+k}(\ZZ/2;\pi_{-k}\bH(M)) \\
&@>>>\pi_{r-1}\wdt{\TOP}^b(M\!\times\!\RR^k)
@>>>\pi_{r-1}\wdt{\TOP}^b(M\!\times\!\RR^{k+1})
@>>>\cdots\endalign
$$
\endproclaim
\msk
The inclusion 
$\bigcup_k\TOP^b(M\times\RR^k)\to \bigcup_k\wdt\TOP^b(M\times\RR^k)$
is a homotopy equivalence \cite{WW1, 1.14}.  Together with 
1.4.2, this shows for example that 
$$\wdt{\TOP}(M)\simeq \bigcup_k\TOP^b(M\times\RR^k)$$
if $M$ is simply connected, because 
then  $\pi_k\bH(M)=0$ for $k\le0$. 
\msk
Another way to relate $\wdt{\TOP}(M)$ and
$\bigcup_k\TOP^b(M\times\RR^k)$  is to use a
filtered version of  Postnikov's method for making highly connected
covers. Let $X$  be a fibrant Space with a filtration by fibrant subSpaces
$$X(0)\quad\subset\quad X(1)\quad\subset\quad X(2)\quad
\subset\dots$$
so that $X$ is the union of the $X(k)$. Call an $i$--simplex in $X$ {\it positive} 
if its characteristic map $\Delta^i\to X$ is 
filtration--preserving, i\.e\., takes the $k$--skeleton of 
$\Delta^i$ to $X(k)$. The positive simplices form a subSpace
${}^{\text{pos}}X$ of $X$, and we let ${}^{\text{pos}}X(k):=
{}^{\text{pos}}X\cap X(k)$. Then ${}^{\text{pos}}X(k)$ is fibrant.
If $X(0)$ is based, then 
$$\pi_i({}^{\text{pos}}X(k),{}^{\text{pos}}X(k-1))
@>{\quad\cong\quad}>>\pi_i(X(k),X(k-1))$$
for $i\ge k$, and $\pi_i({}^{\text{pos}}X(k),{}^{\text{pos}}X(k-1))=*$
for $i<k$.  For example: \linebreak if $X(k)=*$ for $k\le m$ and $X(k)=X$ 
for $k>m$, then ${}^{\text{pos}}X$ is the $m$--connected 
Postnikov cover of $X$. And if $X(k)$ is 
$\TOP^b(M\times\RR^k)/\TOP(M)$, then 
$${}^{\text{pos}}X\,\simeq\wdt\TOP(M)/\TOP(M)\,.$$
See \cite{WW1, 4.10} for a more precise statement, and a proof. 
\msk
\proclaim{1.4.3. Corollary \cite{Ha}} There exists a 
spectral sequence with $E^1$--term given by 
$E^1_{pq}=\pi_{q-1}(\C(M\times I^p))$~, converging 
to the homotopy groups of 
$$\wdt{\TOP}(M)/\TOP(M)\,.$$
\endproclaim
\msk
Hatcher also described $E^2_{pq}$ for $p+n\gg q$.  What he found is 
explained by  the next theorem, which uses 
naturality of the pos--construction and the results of \S1.3. 
\msk
\proclaim{1.4.4. Theorem \cite{WW1, Thm\. C}} There exists a 
homotopy commutative cartesian square of the form 
\msk
$$\CD
\matrix  \wdt{\TOP}(M)/\TOP(M)  \\
\simeq \\
{}^{\text{pos}}\!\left(\bigcup_k\!\TOP^b(M\!\times\!\RR^k)/\TOP(M)\right) 
\endmatrix 
@>>> \Omega\8(\bH^s(M)\ho) \\
@VV\cap V   @VVV  \\
\bigcup_k\TOP(M\!\times\!\RR^k)/\TOP(M) 
@>\Phi>>\Omega\8(\bH(M)\ho) 
\endCD$$
\msk
where $\bH^s(M)$ is the 0--connected cover of $\bH(M)$,  the 
right--hand vertical arrow is induced by the canonical map 
$\bH^s(M)\to \bH(M)$, and $\Phi$ is the map from 1.3.7. 
 \endproclaim
\msk
By 1.3.7,  lower and hence upper horizontal arrow in 1.4.4 are
$j$--connected for $j$ in the topological concordance stable range of 
$M$. 
\bsk
\head 1.5. $h$--structures and $h$--cobordisms \endhead 
Our main goal  in this section is to construct a {\it Whitehead torsion} map
$w:\sS(M)\la
\sH\8(M)$,  and a simple version,
$\sS^s(M)\la(\sH\8)^s(M)=\Omega\8\bH^s(M)$,  which makes the following
diagram homotopy commutative (see 1.1.1):
$$\CD
\wdt{\TOP}(M)/\TOP(M) @>\delta>> \sS^s(M) \\
@V{1.4.4}VV  @VV w V \\
\Omega\8(\bH^s(M)\ho)@>\text{transfer}>>\Omega\8\bH^s(M)\,. 
\endCD\tag1.5.1$$
The map $\delta$ comes from $\wdt{\TOP}(M)/\TOP(M)\to\sS^s(M)\to
\wdt{\sS}^s(M)$, a homotopy fiber sequence. In our description of
$w$, we assume for simplicity that $M^n$ is closed. 
\msk
Let $Z\subset \sS(M)$ be a finitely generated subSpace, that is to
say, $|Z|$ is compact. Let 
$p:E(1)\to |Z|$ be the tautological bundle whose fiber over some vertex $(N,f)$, 
for example, is $N$. (Here $N^n$ is closed and $f:N\to M$ is a homotopy equivalence.) 
Let $E(2)=M\times|Z|$. We have a canonical fiber homotopy
equivalence $\lambda:E(1)\to E(2)$ over $|Z|$. 

Let $\tau_1$ and $\tau_2$ be the vertical tangent bundles of $E(1)$ and $E(2)$,
respectively. Choose $k\gg0$, and a $k$--disk bundle $\xi$ on $E(1)$ with
associated euclidean bundle $\xi\nat$, and an isomorphism $\iota$ of euclidean bundles $\tau_1\oplus\xi\nat\cong\lambda^*(\tau_2\oplus\ep^k)$. Let $E(1)^{\xi}$ be the total 
space of the disk bundle $\xi$; this fibers over $|Z|$. Immersion theory 
\cite{Gau} says that $\lambda$ and $\iota$ together determine up to contractible choice 
a fiberwise codimension zero immersion, over $|Z|$, from $E(1)^{\xi}$ to $E(2)\times\RR^k$. 
We can arrange that the image of this fiberwise immersion is contained in  $E(2)\times\BB^k$ where $\BB^k\subset\RR^k$ is the open unit ball.  Also, by choosing $k$ sufficiently large 
and using general position arguments, we can arrange that the fiberwise 
immersion is a fiberwise embedding. In this situation, the closure of  
$$E(2)\times\DD^k\,\minus\,\im(E(1)^{\xi})$$
is the total space of a fibered family of $h$--cobordisms over 
$|Z|$, with fixed base $M\times\SS^{k-1}$. This family is classified by a map 
$Z\to\sH(M\times\SS^{k-1})$. Letting $k\to\infty$ we have
$$w_Z:Z\,\,\la\,\,\hocolim_k\sH(M\times\SS^{k-1})\,\simeq\,\sH\8(M)\,,$$
a map well defined up to contractible choice. Finally view $Z$ as a variable, 
use the above ideas to make a map $w$ from the homotopy colimit 
of the various $Z$ to $\sH\8(M)$, and note that $\hocolim Z\simeq\sS(M)$. 
This completes the construction. 
\msk
It is evident that $w$ takes $\sS^s(M)$, the Space 
of $s$--structures on $M$, to the Space of
$s$--cobordisms, 
$\Omega\8\bH^s(M)$. Homotopy commutativity  of \thetag{1.5.1} is 
less evident,
but we omit  the proof. 
\msk
\proclaim{1.5.2. Corollary \cite{BuLa2}, \cite{BuFi1}} In the topological concordance stable 
range for $M$, and localized at odd primes, there is a product decomposition 
$$\sS^s(M)\simeq \wdt{\sS}^s(M)\times\wdt{\TOP}(M)/\TOP(M)\,.$$
\endproclaim
\msk
\demo{Proof} The left--hand vertical map in \thetag{1.5.1} is a
homotopy  equivalence in the concordance stable range, and the
lower horizontal map  is a split monomorphism in the homotopy
category, at odd primes. Consequently  the homotopy fiber of 
$$\sS^s(M)\hookrightarrow \wdt{\sS}^s(M)$$
is a retract up to homotopy of $\sS^s(M)$ (localized at odd primes, in the concordance 
stable range). \qed
\enddemo
\msk
\example{1.5.3. Remark} Suppose that $M^n$ is compact with boundary. 
Let $Z\subset\sS(M)$ be finitely generated. A modification
of the construction above gives $w_Z$ from $Z$ to 
$\hocolim_k\sH(M\times\SS^{k-1})$ and then
$w:\sS(M)\to \sH\8(M)$.  

Now let $\sT_n(M)$ be the Space of pairs $(N,f)$ where $N^n$ is a compact 
manifold and $f:N\to M$ is {\it any} homotopy equivalence, not subject to boundary conditions.  
Let $Z\subset \sT_n(M)$ be finitely generated. Another  modification of 
the construction described above gives $w_Z$ from $Z$ to
$\hocolim_k\sH(\partial(M\times\SS^{k-1}))$, and then $w:\sT_n(M)\to\sH\8(M)$ because again
$$\hocolim_k\sH(\partial(M\times\SS^{k-1}))\,\simeq\,\sH\8(M)\,.$$
It follows that $w:\sS(M)\to\sH\8(M)$ extends to 
$w:\sT_n(M)\to\sH\8(M)$.  Let $\G\nat(M)$ 
be the Space of {\it all} homotopy equivalences $M\to M$, not subject to 
boundary conditions. Then $\G(M_0)\simeq \G\nat(M)\hookrightarrow \sT_n(M)$, and we can 
think of $w$ as a map from $\G(M_0)$ to $\sH\8(M)$. This is the map 
promised directly after 1.2.4. 
\endexample
\bsk
\head 1.6. Diffeomorphisms \endhead
Suppose that $M^n$ is a closed topological manifold, $n>4$,
with tangent (micro)bundle $\tau$.  Morlet \cite{Mo1},  \cite{Mo2}, \cite{BuLa3},
\cite{KiSi} proves that the forgetful map from 
a suitably defined Space of smooth structures on $M$, denoted $V(M),$  to a
suitably defined Space $V(\tau)$ of vector bundle structures on $\tau$ is a 
weak homotopy equivalence. Earlier
Hirsch and Mazur \cite{HiMa} had proved that the map in question 
induces a bijection on $\pi_0$. The Space $V(\tau)$ is homotopy 
equivalent to the homotopy fiber of the inclusion 
$\map(M,B\!\Or(n))\to \map(M,B\!\TOP(n))$ over $\hat\tau$. 

The Space of smooth structures on $M$ can be defined as 
the disjoint union of Spaces $\TOP(N)/\DIFF(N)$, where 
$N$ runs through a set of representatives of diffeomorphism 
classes of smooth manifolds homeomorphic to $M$. Therefore 
Morlet's theorem gives, in the case where $M$ is smooth, 
a homotopy fiber sequence 

$$\DIFF(M)@>>>\TOP(M)@>a>>V(M)\,.\tag1.6.1$$

The map $a$ is obtained from an action of $\TOP(M)$ on $V(M)$ 
by evaluating at the base point of $V(M)$. 
The homotopy fiber sequence  \thetag{1.6.1} remains meaningful when $M$ 
is smooth compact with boundary;  in this case allow only 
vector bundle structures on $\tau$ extending the standard 
structure over $\partial M$. 
\msk
Traditionally, 1.6.1 has been an excuse for neglecting $\DIFF(M)$ 
in favor of $\TOP(M)$. But concordance theory has changed that. 
See \S3.  It is therefore best to develop 
a theory of smooth automorphisms parallel to the theory of 
topological automorphisms where possible. 
For example, 1.2.1--3, 1.3.1--2,  1.3.5--8, 1.4.2--4 and 1.5.1--3 
have smooth analogs; but 1.2.4 does not. 

\noindent {\bf Notation:} A subscript $d$ will often be used  to 
indicate smoothness, as in $\sH_d(M)$ for a Space of smooth
$h$--cobordisms (when $M$ is smooth).

\bsk
\head  {\bf 2. $L$--theory and structure Spaces} \endhead
The major theorems in this chapter are to be found in 
\S2.3 and \S2.5.  Sections 2.1, 2.2 and 2.4 introduce concepts 
needed to state those theorems. 
\msk 
\head 2.1. Assembly \endhead
\definition{2.1.1. Definitions} Fix a space $Y$. Let $\sW_Y$ be the category 
of spaces over $Y$. A morphism in $\sW_Y$ is a 
{\it weak homotopy equivalence} if the underlying map of spaces is a 
weak homotopy equivalence. A commutative square in $\sW_Y$ is 
{\it cartesian} if the underlying square of spaces is cartesian.

A functor $\bJ$ from $\sW_Y$ to CW--spectra \cite{A, III} is 
{\it homotopy invariant} if it takes weak homotopy equivalences 
to weak homotopy equivalences. It is {\it excisive} if, in addition,
it takes cartesian squares to cartesian squares, takes $\emptyset$ 
to a contractible spectrum, and satisfies a wedge axiom, 
$$\vee_i \bJ(X_i\to Y) @>\simeq>> \bJ(\amalg_iX_i\to Y)\,.$$
\enddefinition
\msk
\proclaim{2.1.2. Proposition} For every homotopy invariant 
functor $\bJ$ from $\sW_Y$ to CW--spectra, there exist an excisive
functor $\bJ\upr$ from $\sW_Y$ to CW--spectra and a natural 
transformation $\alpha_{\bJ}:\bJ\upr\to\bJ$ such that 
$\alpha_{\bJ}:\bJ\upr(X)\to \bJ(X)$
is a homotopy equivalence whenever $X$ is a point (over $Y$). 
\endproclaim
\msk
The natural transformation $\alpha_{\bJ}$ is essentially characterized by these 
properties;  it is called the {\it assembly}. It is the best approximation
(from the left) of $\bJ$ by an excisive functor.  We write 
$\bJ\lpr(X)$ for the homotopy fiber of $\alpha_{\bJ}:\bJ\upr(X)\to
\bJ(X)$, for any $X$ in $\sW_Y$. 
\msk
\remark{Remark} Assume that $X$ and $Y$ above are homotopy 
equivalent to CW--spaces. If $Y\simeq*$, then 
$\bJ\upr(X)\simeq X_+\wedge \bJ(*)$ by a chain of natural 
homotopy equivalences. If $Y\not\simeq*$, build a quasi--fibered
spectrum on $Y$  with fiber $\bJ(y\hookrightarrow Y)$ over $y\in Y$. 
Pull it back to $X$ using $X\to Y$. Collapse the zero section, 
a copy of $X$. The result is $\simeq\bJ\upr(X)$.
\endremark
\msk
The assembly concept is due to Quinn, \cite{Qun1},  \cite{Qun2}, \cite{Qun3}. See also 
\cite{Lo}. For a proof 
of 2.1.2, see \cite{WWa}. For applications to block 
$s$--structure Spaces we need the case where
$Y=\RR P\8$ and $\bJ=\bL^s\bul$ is the $L$--theory functor
$X\mapsto\bL\bul^s(X)$ which associates to $X$ the 
$L$--theory spectrum of $\ZZ\pi_1(X)$, say in the 
description of Ranicki \cite{Ra2 \S13}. Note the following technical points.
\roster
\itb $\bL^s\bul$ is the quadratic $L$--theory 
with decoration $s$ which we make 0--connected (by force). 
\itb Because of 2.1.1 we have no use for a base point in $X$.
This makes it harder to say what  
$\ZZ\pi_1(X)$ should mean. For details see \cite{WWa, 3.1}, where 
$\ZZ\pi_1(X)$ or more precisely $\ZZ^w\pi_1(X)$ 
is a ringoid with involution depending on the double cover 
$w:X\td\to X$ induced by $X\to\RR P\8$.
\endroster
\bsk
\head 2.2. Tangential invariants  \endhead 
Geometric topology tradition requires that any 
classification of you--name--it structures on a manifold or 
Poincar\'e space \cite{Kl} be
accompanied by a classification of analogous structures on the 
normal bundle or Spivak normal fibration \cite{Spi}, \cite{Ra3}, \cite{Br2} of the manifold or 
Poincar\'e space.  We endorse this.  However, we find tangent bundle
language more convenient than normal bundle language. The
constructions  here in \S2.2 will also be used in \S3 and \S4. 
\remark{Terminology} When we speak of a {\it stable fiber homotopy
equivalence} between euclidean bundles $\beta$ and $\gamma$ on 
a space $X$, we mean a fiber homotopy equivalence over $X$ between the 
spherical fibrations associated with $\beta\oplus\ep^k$ and
$\gamma\oplus\ep^k$~, respectively, for some $k$. 
\endremark
\msk
Let $M^n$ be a closed topological manifold  with a
choice \cite{Kis}, \cite{Maz2} of euclidean tangent bundle $\tau$. 
An {\it $h$--structure} on $\tau$ is a pair $(\xi,\phi)$ 
where $\xi^n$ is a euclidean bundle on $M$ and $\phi$ is a stable 
fiber homotopy equivalence from (the spherical 
fibration associated with) $\xi$ to $\tau$. The $h$--structures on $\tau$
and  their isomorphisms form a groupoid. Enlarge the groupoid to a 
simplicial groupoid by allowing families parametrized by $\Delta^k$~,
and let 
$$\sS(\tau):=\,\,\text{diagonal nerve of the simplicial groupoid}. $$
Then $\sS(\tau)\,\,\simeq\,\,
\hofiber\,[\,\map(M,B\!\TOP(n))\hookrightarrow 
\map(M,B\!\G)\,]$~, 
where $\hat\tau$ serves as base point in
$\map(M,B\!\TOP(n))$ and $\map(M,B\!\G)$.  There is a 
{\it tangential invariant map}, well defined up to contractible choice,
$$\nabla:\sS(M)\la \sS(\tau)\,.$$
Sketchy description: A homotopy equivalence $f:N\to M$ determines 
up to contractible choice a stable fiber homotopy equivalence 
$\psi$ from $f^*\nu(M)$ to $\nu(N)$ because Spivak normal fibrations
\cite{Spi}, 
\cite{Ra3}, \cite{Br3},  \cite{Wa5} are homotopy invariants. It then
determines up to  contractible choice a stable fiber homotopy equivalence
$\psi^{\text{ad}}:\tau(N)\to f^*\tau(M)$. Now choose a euclidean
bundle $\xi^n$ on $M$ 
and a stable fiber homotopy equivalence $\phi:\xi\to \tau(M)$ 
{\it together} with an isomorphism $j:f^*\xi\to\tau(N)$ and a homotopy 
$f^*\phi\simeq \psi^{\text{ad}}\cdot j$. This is a contractible 
choice. Let $\nabla$ take $(N,f)$ to $(\xi,\phi)$. (This defines 
$\nabla$ on the 0--skeleton; using the same ideas, complete 
the construction of $\nabla$ by induction over skeletons.)
\msk
When $\partial M\ne\emptyset$, define $\sS(\tau)$
in such a way that it is homotopy equivalent to the homotopy 
fiber of 
$$\map_{\text{rel}}(M,B\!\TOP(n))@>\subset>>\map_{\text{rel}}(M,B\!\G)$$
where $\map_{\text{rel}}$ indicates maps from $M$ which on $\partial M$ 
agree with $\hat\tau$. Again  
$\hat\tau$ serves as base point everywhere. 
\msk
When $\partial M\ne\emptyset$
 and $\partial_+M\subset\partial M$ is specified, with tangent 
bundle $\tau'$ of fiber dimension $n-1$, we define $\sS(\tau,\tau')$
in such a way that it is homotopy equivalent to the homotopy 
fiber of 
$$\map_{\text{rel}}((M,\partial_+M)\,,\,(B\!\TOP(n),B\!\TOP(n-1))
@>\subset>>\map_{\text{rel}}(M,B\!\G)$$
where $\map_{\text{rel}}$ indicates maps from $M$ which on $\partial_-M$ 
agree with the classifying map for the tangent bundle of 
$\partial_-M$.  The pair $(\hat\tau,\hat\tau')$
serves as base point. There is a {\it tangential invariant map}
$\nabla:\sS(M,\partial_+M)\to \sS(\tau,\tau')$ which fits 
into a homotopy commutative diagram where the rows are 
homotopy fiber sequences:
$$\CD
\sS(M)@>\subset>>\sS(M,\partial_+M)@>>> \sS(\partial_+M) \\
@VV\nabla V   @VV\nabla V   @VV\nabla V  \\
\sS(\tau)@>\subset>>\sS(\tau,\tau')@>>>\sS(\tau')\,.
\endCD$$
\msk
\example{2.2.1. Illustration} Suppose that $(M,\partial_+M)=(N\times I,N\times 1)$ 
where $N^{n-1}$ is compact. Then $\sS(\tau,\tau')$ is 
homotopy equivalent to the homotopy fiber of 
$$\map_{\text{rel}}(N,B\!\TOP(n-1))\hookrightarrow\map_{\text{rel}}(N,B\!\TOP(n))\,.$$
\endexample
\msk
Finally we need a space $\wdt{\sS}(\tau)$ of {\it stable} 
$h$--structures on $\tau$ (assuming again that $\tau$ is the 
tangent bundle of a compact $M$). Define this as the space 
of pairs $(\xi,\phi)$ where $\xi^p$ is a euclidean bundle on $M$, 
of arbitrary fiber dimension $p$, and $\phi$ is a stable fiber 
homotopy equivalence $\xi\to\tau$, represented by an actual 
fiber homotopy equivalence between the spherical fibrations 
associated with $\xi\oplus\ep^{k-p}$ and $\tau\oplus\ep^{k-n}$
for some large $k$. Then 
$$\wdt{\sS}(\tau)\,\,\simeq\,\,\hofiber[\map_{\text{rel}}(M,B\!\TOP)\to\map_{\text{rel}}(M,B\!\G)]
$$
which is homotopy equivalent to the space of based 
maps from $M/\partial M$ to $\G\!/\!\TOP$. 
Again there is a tangential invariant map--- 
better known, and in this case more easily described, as 
the {\it normal invariant map} :
$$\nabla:\wdt{\sS}^s(M)\to \wdt{\sS}(\tau)\,.$$
\bsk
\head 2.3. Block $s$--structures and $L$--theory  \endhead
\proclaim{2.3.1.  Fundamental Theorem of Surgery 
{\rm (Browder, Novikov, Sullivan, Wall, Quinn, Ranicki)}}
For compact $M^n$ with tangent bundle $\tau$, where $n>4$,
there exists a homotopy commutative square of the 
form 
$$
\CD
\wdt{\sS}^s(M)@>\nabla>> \wdt{\sS}(\tau) \\
@VV\simeq V   @VV\simeq V \\
\Omega^{\infty+n}((\bL\bul^s)\lpr(M))@>\text{forget}>> 
\Omega^{\infty+n}((\bL\bul^s)\upr(M))\,. \\
\endCD$$
\endproclaim
\msk
References: 
\cite{Br3}, \cite{Br4}, \cite{Nov} for the smooth analog, \cite{Rou1},
\cite{Su1}, \cite{Su2}, \cite{Wa1,\S10}, \cite{ABK} for the PL case, 
and \cite{KiSi} for the topological case, all without explicit use of 
assembly~;  \cite{Ra4}, \cite{Ra2}, 
\cite{Qun1} for formulations with  assembly. 
\msk
{\it Illustration.} 2.3.1 gives 
$\wdt{\sS}^s(\SS^n)\,\simeq\,\Omega\8\bL\bul^s(*)\simeq \G\!/\TOP$
for $n>4$.
(Remember that $\bL\bul^s(*)$ is 0--connected by definition, \S2.1.)
The honest structure space is 
$$\sS(\SS^n)\,\simeq\,\G(\SS^n)/\TOP(\SS^n)$$
(this uses the Poincar\'e conjecture). The inclusion 
$\sS(\SS^n)\to\wdt{\sS}^s(\SS^n)$ becomes the inclusion 
of $\G(\SS^n)/\TOP(\SS^n)$ in 
$\hocolim_k\G(\SS^k)/\TOP(\SS^k)\,\simeq\,\G\!/\TOP$~; in particular, 
it is $n$--connected. 
\bsk
\head 2.4. $\SS^1$--stabilization \endhead
$\SS^1$--stabilization is a method for making new homotopy 
functors out of old ones. It was introduced in \cite{Ra5} 
and applied in \cite{AnPe}, \cite{HaMa} in a special case (the one we will 
need in \S2.6 just below). It is motivated by the definition of the negative
$K$--groups in \cite{Ba}.  

Let $\bJ$ be  a homotopy functor from $\sW_Y$ (see 2.1.2) to 
CW--spectra.  Let $\SS^1(+)$ and $\SS^1(-)$ be upper half 
and lower half of $\SS^1$, respectively. For $X$ in $\sW_Y$ let 
$\sigma\bJ(X)$ be the homotopy pullback of 
$$\dsize\frac{\bJ(X\times\SS^1(+))}{\bJ(X\times *)}\la 
\dsize\frac{\bJ(X\times\SS^1)}{\bJ(X\times *)} \longleftarrow
\dsize\frac{\bJ(X\times\SS^1(-))}{\bJ(X\times *)}\,.$$
Note that 
$\sigma\bJ(X)\,\simeq\, \Omega[\,\bJ(X\times\SS^1)/\bJ(X\times *)\,]$,
and that $\sigma\bJ$ is a homotopy functor from $\sW_Y$ to
CW--spectra. There are natural transformations 
$$\bJ(X)\la\frac{\bJ(X\times\SS^0)}{\bJ(X\times *)}
\la\sigma\bJ(X)\,,$$
the first induced by the inclusion 
$x\mapsto (x,-1)$ of $X$ in $X\times\SS^0$,  and the second induced 
by the inclusions of $\SS^0$ in $\SS^1(-)$ and $\SS^1(+)$.
Let $\psi:\bJ(X)\to \sigma\bJ(X)$ be the composition. Finally let
$\sigma\8\bJ(X)$ be the homotopy colimit of the 
$\sigma^k\bJ(X)$ for $k\ge0$, using the maps
$\psi:\sigma^{k-1}\bJ(X)\to\sigma(\sigma^{k-1}\bJ(X))$ to 
stabilize. We
call $\sigma\8\bJ$ the  {\it $\SS^1$--stabilization} of $\bJ$.
\msk
\head 2.5.  Bounded $h$--structures and $L$--theory  \endhead
Let $M^n$ be compact. A {\it bounded} $h$--structure
on $M\times\RR^k$ is a pair $(N,f)$ where $N^{n+k}$ is a manifold 
and $f:N\to M\times\RR^k$ is a {\it bounded} homotopy equivalence
restricting to a homeomorphism $\partial N\to\partial M\times\RR^k$.
(That is, there exist $c>0$ and $g:M\times\RR^k\to N$ and homotopies
$h:fg\simeq\id$, $j:gf\simeq\id$ such that the sets 
$\{p_2h_t(x)\mid t\in I\}$, $\{p_2fj_t(y)\mid t\in I\}$ have diameter
$<c$ for all $x\in M\times\RR^k$ and $y\in N$; moreover 
$h_t$, $j_t$ agree with the identity maps on $\partial M\times\RR^k$ 
and $\partial N$, respectively.) References: \cite{AnPe}, \cite{FePe}. 

A  Space $\sS^b(M\times\RR^k)$ of such 
bounded $h$--structures can be constructed in the 
usual way, as the diagonal nerve of a simplicial groupoid. 
There is a homotopy fiber sequence
$$\TOP^b(M\times\RR^k)\to \G^b(M\times\RR^k)\to 
\sS^b(M\times\RR^k)$$
where $G^b(M\times\RR^k)$ is the Space of bounded homotopy 
automorphisms of $M\times\RR^k$, relative to $\partial M\times\RR^k$. 

Again we need a {\it tangential invariant map} 
$\nabla:\sS^b(M\times\RR^k)\to
\sS(\tau\times\ep^k)$
where $\tau=\tau(M)$ and $\tau\times\ep^k$ is the 
tangent bundle of $M\times\RR^k$.  
Its definition resembles that of the tangential
invariant maps in \S2.2. Additional subtlety:
one needs to know that {\it open} Poincar\'e spaces and open 
Poincar\'e pairs  have Spivak normal fibrations which are invariants of 
proper homotopy type. See \cite{Tay}, \cite{Mau}, \cite{FePe}, \cite{PeRa}. Note
$$\sS(\tau\times\ep^k)\,\,\simeq\,\,\,\hofiber\,
[\map_{\text{rel}}(M,B\!\TOP(n+k))\to \map_{\text{rel}}(M,B\!\G)\,]$$
where $\map_{\text{rel}}$ indicates maps which on $\partial M$ 
agree with $\hat\tau$. 
\msk
Let $\bL\bul^{\langle -\infty\rangle}(X)$ be the 0--connected cover 
of $(\sigma\8\bL\bul^s)(X)$, in the notation of \S2.4. Here $X$ is a space
over $\RR P\8$.
\msk
\proclaim{2.5.1. Theorem}
For compact $M^n$ with tangent bundle $\tau$,  where $n>4$,
there exists a homotopy commutative square of the 
form 
$$
\CD
\bigcup_k\sS^b(M\times\RR^k)@>\nabla>> \bigcup_k\sS(\tau\times\ep^k) \\
@VV\simeq V   @VV\simeq V \\
\Omega^{\infty+n}((\bL\bul^{\langle-\infty\rangle})\lpr(M))@>\text{forget}>> 
\Omega^{\infty+n}((\bL\bul^{\langle-\infty\rangle})\upr(M))\,. \\
\endCD$$
\endproclaim
\msk
\remark{Remark} 
$(\bL\bul^{\langle-\infty\rangle})\upr\simeq(\bL\bul^s)\upr$.
\endremark
\bsk
\head {\bf 3. Algebraic K-theory and structure Spaces} \endhead
\head 3.1.  Algebraic K-theory of Spaces  \endhead
Waldhausen's homotopy functor $\bA$ from spaces to CW--spectra 
is a composition $\bK\cdot\sR$, where 
$\sR$ is a functor from spaces to {\it categories with 
cofibrations and weak equivalences} alias {\it Waldhausen categories},
and $\bK$ is a functor from Waldhausen categories to CW--spectra. 

For a space $X$, let $\sR(X)$ be the Waldhausen category of 
homotopy finite retractive spaces over $X$.  The objects of $\sR(X)$ are 
spaces $Z$ equipped with maps
$$\qquad\qquad\qquad
Z \overset r\to{\overset{\dsize\la}\to{\underset i\to\longleftarrow}} X
\qquad\qquad\qquad(ri=\id_X)$$ 
subject to a finiteness condition. Namely,  $Z$ must be 
homotopy equivalent, relative to $X$, to a relative CW--space built 
from $X$ by attaching a finite number of cells. 
A morphism in $\sR(X)$ is a map relative to and over $X$.  We call 
it a weak equivalence if it is a homotopy equivalence relative to $X$,
and a cofibration if it has the homotopy extension property 
relative to $X$.  See \cite{Wah3, ch.2} for more information. 
\msk
In \cite{Wah3}, \cite{Wah1}, Waldhausen 
associates with any Waldhausen category $\sC$ a connective spectrum 
$\bK(\sC)$, generalizing Quillen's construction \cite{Qui} of 
the $K$--theory spectrum of an exact category. For us it is 
important that $\bK(\sC)$ comes with a  map reminiscent 
of ``group completion'',
$$|w\sC|\hookrightarrow \Omega\8\bK(\sC)$$
where $w\sC$ is the category of weak equivalences in $\sC$ and 
$|w\sC|$ is its classifying space (geometric realization of the nerve). 
\bsk
\head 3.2. Algebraic $K$--theory of spaces, and 
$h$--cobordisms
\endhead 
Let $M^n$ be compact, $n\ge 5$, with fundamental
group(oid) $\pi$.  The $s$--cobordism theorem due to Smale
\cite{Sm}, \cite{Mi2}   in the simply connected smooth case and
Barden--Mazur--Stallings \cite{Ke} in the nonsimply connected
smooth case states that $\pi_0\sH_d(M)$ is isomorphic to the
Whitehead group of $\pi$, that is,
$K_1(\ZZ\pi)/\{\pm\pi^{\text{ab}}\}$. See
\cite{RoSa} for the PL version and \cite{KiSi, Essay III} for the $\TOP$
version. Cerf \cite{Ce} showed that $\pi_1\sH_d(M)$ is trivial when $M$
is smooth, simply connected and $n\ge5$, and Rourke \cite{Rou2}
established the analogous statement in the PL category. 
In the early 70's Hatcher and Wagoner \cite{HaWa}, working with a
smooth but possibly  nonsimply connected $M$,  constructed a surjective homomorphism
from $\pi_1\sH_d(M)$  to a certain quotient of $K_2(\ZZ\pi)$, and they were able to
describe the kernel of that homomorphism in terms of $\pi$ and $\pi_2(M)$. See also
\cite{DIg}.  These results follow from Waldhausen's
theorem 3.2.1, 3.2.2 below, which describes the homotopy 
types of $\sH(M)$ and $\sH_d(M)$ in a stable range, in algebraic
$K$--theory terms.  The size of the stable range is estimated by 
Igusa's stability theorem,  1.3.4.
\remark{Remark} Note that the block $s$--cobordism
Space $\wdt{\sS}^s(M\times I,M\times 1)$ is not a very useful 
approximation to $\sH^s(M)$~, because it is contractible (either by a 
relative version of 2.3.1 which we did not state, or by a direct 
geometric argument). 
Hence surgery theory as in \S2 does not elucidate the homotopy 
type of $\sH^s(M)$. 
\endremark
\msk
For compact $M^n$ write 
$\sH(\tau(M)):=\sS(\tau(M\times I),\tau(M\times 1))$ so that 
there is a tangential invariant map $\nabla:\sH(M)\to\sH(\tau(M))$.
See 2.2.1. There is a stabilization map from
$\sH(\tau(M))$ to $\sH(\tau(M\times I))$, analogous to 
the stabilization map from $\sH(M)$ to $\sH(M\times I)$. We let
$\tau=\tau(M)$  and $\sH^{\infty}(\tau)=\hocolim_k\sH(\tau(M\times I^k))$
and obtain, since $\nabla$ commutes with stabilization, 
$$\nabla:\sH^{\infty}(M)\to \sH^{\infty}(\tau)\,.$$ 
The following result is essentially 
contained in \cite{Wah2}.
\msk
\proclaim{3.2.1. Theorem {\rm(Waldhausen)}} There exists a
homotopy  commutative square
$$
\CD
\sH^{\infty}(M) @>\nabla>>\sH^{\infty}(\tau) \\
@VV\simeq V    @VV\simeq V \\
\Omega\8(\bA\lpr(M))@>\text{forget}>>\Omega\8(\bA\upr(M))\,.
\endCD
$$
\endproclaim
\msk
\example{3.2.2. Remark} Suppose that $M$ is smooth. Then
$\sH^{\infty}(M)$ and $\sH^{\infty}(\tau)$ have smooth 
analogues $\sH_d^{\infty}(M)$ and $\sH_d^{\infty}(\tau)$, and 
by smoothing theory there is a homotopy commutative cartesian
square 
\msk
$$\CD
\sH_d^{\infty}(M)@>\nabla>>\sH_d^{\infty}(\tau) \\
@VVV @VVV \\
\sH^{\infty}(M)@>\nabla>>\sH^{\infty}(\tau)
\endCD$$
\msk
with forgetful vertical arrows. One shows by direct geometric arguments  
that $\sH_d^{\infty}(\tau)\simeq\Omega\8\Sigma\8(M_+)$ and 
that $\nabla:\sH_d^{\infty}(M)@>>>\sH_d^{\infty}(\tau)$ is 
nullhomotopic. In this way 3.2.1 implies 
$$\sH_d^{\infty}(M)\times \Omega^{\infty+1}\Sigma\8(M_+)
\quad\simeq\quad \Omega^{\infty+1}(\bA(M))\tag3.2.3$$
which is better known than 3.2.1.  Conversely, 3.2.1 can be 
deduced from \thetag{3.2.3} with functor calculus arguments, 
if we add the information 
that  \thetag{3.2.3} comes from a spectrum level splitting, 
$\bH^h_d(M)\vee\Omega\Sigma\8(M_+)
\simeq \Omega\bA(M)$ or equivalently $\Whd(M)\vee\Sigma\8(M_+)\simeq \bA(M)$, where 
$\Whd(M)$ is the delooping of $\bH^h_d(M)$ and $\bH^h_d(M)$ is the 
$(-1)$--connected cover of $\bH_d(M)$. 
We prefer formulation 3.2.1 because of its amazing similarity 
with 2.3.1 and 2.5.1. 
\endexample
\msk
References: 3.2.3 is stated in  \cite{Wah2}.  It is reduced in 
\cite{Wah4} to the spectrum level analog of the left--hand 
column in \thetag{3.2.1}, 
$$\bH^h(M)\simeq\bA\lpr(M)~,\tag3.2.4$$
where $\bH^h(M)$ denotes the $(-1)$--connected cover of
$\bH(M)$.) For the proof of \thetag{3.2.4}, see
\cite{Wah3, \S3} and the preprints \cite{WaVo1} and \cite{WaVo2}.
The papers \cite{Stb} and \cite{Cha5} contain results 
closely related to \cite{WaVo1} and \cite{WaVo2}, 
respectively.  A very rough but helpful guide to this vast circle of
ideas is \cite{Wah5}. See also \cite{DWWc}. 
\bsk
\head  {\bf 4. Mixing $L$--theory and algebraic $K$--theory
of spaces}
\endhead
\remark{Introduction} We now have a large amount of indirect knowledge about the 
$s$--structure Space $\sS^s(M)$ for a compact $M$. Namely, 
from the definitions there is a homotopy fiber sequence
$$\wdt{\TOP}(M)/\TOP(M)\la
\sS^s(M)\la \wdt{\sS}^s(M).$$
In 2.3.1 we have an expression for $\wdt{\sS}^s(M)$ in terms of $L$--theory.  
In the concordance stable range, we also have the expression 1.3.7 for 
$$\wdt{\TOP}(M)/\TOP(M)$$
in terms of stabilized concordance theory.  But 3.2.1 expresses 
stabilized concordance theory through the algebraic $K$--theory 
of spaces.
Therefore, in the concordance stable range, $\sS^s(M)$ must be a
concoction of
$L$--theory and algebraic $K$--theory of spaces. It remains to find
out what concoction exactly.  This problem was previously addressed
by  Hsiang--Sharpe  (roughly speaking, using only the Postnikov
2--coskeleton of the algebraic $K$--theory of spaces),  by
Burghelea--Fiedorowicz (rationally),  by Burghelea--Lashof (at odd
primes), by Fiedorowicz--Schw\"anzl--Vogt (at odd  primes); see
references
\cite{HsiSha}, \cite{BuFi1}, \cite{BuFi2}, 
\cite{FiSVo1}, \cite{FiSVo2}, \cite{FiSVo3}, \cite{BuLa2}. In addition, the 
literature contains many results about $\sS^s(M)$ or $\sS(M)$, or the 
differentiable analogs, for specific $M$; see \S6 for a selection 
and further references. 
\msk
Our analysis, Thm\. 4.2.1 and remark 4.2.3 below, is based on the
following idea.  Using 3.2.1, Poincar\'e duality notions, and a more  
algebraic description of $w$ in \thetag{1.5.1}, we
find that $w$ lifts to an {\it equivariant Whitehead torsion} map 
$$w\sha:\sS^s(M)\la \Omega\8(\bH^s(M)\hf) \simeq
\Omega\8((\bA^s\lpr(M)\hf))$$ where $(\text{---})\hf$ indicates
homotopy fixed points for a certain  action of $\ZZ/2$. This
refinement of $w$ fits into a homotopy commutative diagram
whose top portion refines \thetag{1.5.1}, 
$$\CD
\wdt{\TOP}(M)/\TOP(M)  @>>> \Omega\8((\bA^s\lpr(M)\ho)) \\
@VVV   @VV\text{norm} V  \\
\sS^s(M)@>w\sha>> \Omega\8((\bA^s\lpr(M)\hf)) \\
@VVV @VVV \\
\wdt{\sS}^s(M) @>>> \Omega\8((\bA^s\lpr(M)\thf))\,.
\endCD$$
The right--hand column is a homotopy fiber sequence 
of infinite loop spaces which we will say 
more about below.  The left--hand column is also 
a homotopy fiber sequence.  The
upper horizontal arrow is highly connected. Hence the 
lower square of the diagram is {\it approximately} 
cartesian.  Since we have an algebraic description for the 
lower square with $\sS^s(M)$ deleted, 
we obtain an approximate algebraic description of $\sS^s(M)$. 

Curiously, the map $w\sha$ does not have an
easy analog in the smooth category.  
See however \S4.3. 
\msk
This work is still in progress. Currently available:
\cite{WW1}, \cite{WW2}, \cite{WWa}, \cite{WWx}, \cite{WWd}, \cite{WWp}, \cite{We1}.
The papers \cite{DWW}, \cite{DWWc} are closely related and use identical
technology. 
\endremark
\bsk
\head 4.1. $LA$--theory \endhead
We describe a functor $\bL\bA^h\bul$ from $\sW_{B\!\G}\times\NN$ to 
CW--spectra. Here $\sW_{B\!\G}$ is the category of spaces over $B\!\G$ 
(alternatively, spaces equipped with a stable spherical fibration) and 
$\NN$ is regarded as a category with exactly one morphism 
$m\to n$ if $m\le n$, and no morphism $m\to n$ if $m>n$. 
For fixed $n$, the functor $\bL\bA^h\bul(\text{---},n)$
is a homotopy functor. It is a composition $\bF_n\cdot\sR\8$, 
where $\sR\8$ is a functor from $\sW_{B\!\G}$ to the 
category of Waldhausen categories with Spanier--Whitehead (SW)
product \cite{WWd}, and $\bF_n$ is a functor from certain 
Waldhausen categories with SW product to CW--spectra. 

Let $X$ be a space over $B\!\G$. In the Waldhausen category
$\sR(X)$ of 3.1, we have notions of mapping cylinder, mapping cone and 
suspension $\Sigma_X$. Let $\sR\8(X)$ be the colimit of the direct 
system of categories
$$\sR(X)@>\Sigma_X>>\sR(X)@>\Sigma_X>>\sR(X)@>\Sigma_X>>\dots\,.$$
Again, $\sR\8(X)$ is a Waldhausen category. 
We will need additional structure on $\sR\8(X)$ in the shape of an 
{\it SW product} which depends on the reference map from $X$
to $B\!\G$.  The SW product is a 
functor $(Z_1,Z_2)\mapsto Z_1\odot Z_2$ from $\sR\8(X)\times\sR\8(X)$ to 
based spaces. Its main properties are:
\roster
\itb  {\it symmetry}, that is, 
$Z_1\odot Z_2\cong Z_2\odot Z_1$ by a natural involutory 
homeomorphism;
\itb {\it bilinearity}, that is,  for fixed $Z_2$ the functor 
$Z_1\mapsto Z_1\odot Z_2$ takes the zero object to a contractible space
and takes pushout squares where the horizontal arrows are 
cofibrations to cartesian squares;
\itb {\it $w$--invariance}, that is, a weak equivalence 
$Z_1\to Z_1'$ induces a weak homotopy equivalence 
$Z_1\odot Z_2\to Z_1'\odot Z_2$ for any $Z_2$.
\endroster
Modulo technicalities, the definition of $Z_1\odot Z_2$ for 
$Z_1,Z_2$ in $\sR(X)\subset\sR\8(X)$ is as follows. Let $\gamma$ be 
the spherical fibration on $X$ pulled back from $B\!\G$; we can 
assume that it comes with a distinguished section and has 
fibers $\simeq\SS^k$ for some $k$. 
Convert the retraction maps $Z_1\to X$ and $Z_2\to X$ into 
fibrations, with total spaces $Z_1\td$ and $Z_2\td$~; form the 
fiberwise smash product (over $X$) of $Z_1\td$, $Z_2\td$ and 
the total space of $\gamma$; collapse the zero section to a point;
finally apply $\Omega^{\infty+k}\Sigma\8$. 
\msk
Imitating \cite{SpaW},  \cite{Spa,\S8 ex\.F} we say that $\eta\in Z'\odot Z$ 
is a {\it duality} if it has certain nondegeneracy properties. 
For every $Z$ in $\sR(X)$ there exists a $Z'$ and $\eta\in Z'\odot Z$ 
which is a duality; the pair $(Z',\eta)$ is 
determined up to contractible choice by $Z$, and we can say 
that $Z'$ is the {\it dual} of $Z$.  Modulo technicalities, an 
involution on the $K$--theory spectrum  $\bK(\sR\8(X))$ 
results, induced by $Z\mapsto Z'$.  The inclusion 
$\sR(X)\to\sR\8(X)$ induces a homotopy equivalence 
of the $K$--theory spectra, so that we are talking about an 
involution on $\bA(X)$.  For all details,
we refer to \cite{WWd}. The involution on $\bA(X)$ was
first constructed in \cite{Vo}. See also \cite{KVWW2}.   
\msk
More generally, suppose  that $\sB$ is any Waldhausen category with SW
product, satisfying the axioms of \cite{WWd,\S2} which assure 
existence and essential uniqueness of SW duals. Then the $K$--theory 
spectrum $\bK(\sB)$ has a preferred involution. In this setup it is 
also possible to define spectra $\bL\bul(\sB)$ (quadratic $L$--theory),
$\bL\ubul(\sB)$ (symmetric $L$--theory), a forgetful map 
from quadratic to symmetric $L$--theory, and a map 
$$\Xi:\bL\ubul(\sB)\la \bK(\sB)\thf$$
where $\bK(\sB)\thf$ is the mapping cone of the {\it norm map}
\cite{AdCD}, \cite{GrMa},  \linebreak \cite{WW2, 2.4},
$$\bK(\sB)\ho\la \bK(\sB)\hf\,,$$
from the homotopy orbit
spectrum to the homotopy fixed point spectrum of the involution 
on $\bK(\sB)$. The norm map refines the transfer map from 
$\bK(\sB)\ho$ to $\bK(\sB)$. 
\msk
Our constructions
of $\bL\bul(\sB)$ and $\bL\ubul(\sB)$ are bordism--theoretic 
and follow \cite{Ra2} very closely, except that with a view to 
the applications here we need 0--connected versions. 
In particular, if $\sB=\sR\8(X)$, 
then $\bL\bul(\sB)$ is homotopy equivalent to the 0--connected 
cover of the quadratic $L$--theory spectrum (decoration $h$) of the 
ring(oid) with involution $\ZZ\pi_1(X)$. 
Regarding $\Xi$, we offer the following explanations. Let $\sB^*[i]$ be
the category of covariant functors from the poset of nonempty faces
of $\Delta^i$ to $\sB$. Then $\sB^*[i]$ inherits from $\sB$ the 
structure of a Waldhausen category with SW product, with weak 
equivalences and cofibrations defined coordinatewise.  The 
axioms of \cite{WWd,\S2} are still satisfied. Modulo technicalities,  there is a duality involution on $|\,w\sB^*[i]\,|$ for each $i\ge0$, and 
an inclusion of simplicial spaces:
\msk
$$\CD
\qquad\qquad\quad i\quad\mapsto\quad |\,w\sB^*[i]\,|^{h\ZZ/2} \\
@VVV  \\
\qquad\qquad\quad i\quad\mapsto\quad K(\sB^*[i])^{h\ZZ/2}\,.
\endCD$$
The geometric realizations of these simplicial spaces turn out to be 
$\Omega\8$ of $\bL\ubul(\sB)$ and $\bK(\sB)\thf$ respectively.
(Recognition is easy in the first case, harder in the second case.)
The inclusion map of geometric realizations is $\Omega\8$ of $\Xi$~, by 
definition.
\msk
We come to the description of $\bF_n(\sB)$, promised at the 
beginning of this section.
Let $\SS^n_!=\RR^n\cup \infty$ with the involution $z\mapsto
-z$ for $z\in\RR^n$. This has fixed point set 
$\{0,\infty\}\cong\SS^0$.  Let
$\bK(\sB,n):=\SS^n_!\wedge\bK(\sB)$ with  the diagonal
involution. The inclusion of $\bK(\sB)\cong 
\bK(\sB,0)$ in $\bK(\sB,n)$ induces a homotopy 
equivalence of Tate spectra, 
$$\bK(\sB)\thf@>\quad t^n\quad>> \bK(\sB,n)\thf$$
(proof by induction on $n$). 
Write ? for the forgetful map from quadratic $L$--theory 
to symmetric $L$--theory.
Let $\bF_n(\sB)$ be the homotopy pullback of 
\msk
$$\CD
 @. \bK(\sB,n)\hf \\
@.        @VV\cap V \\
\bL\bul(\sB)@>{\quad t^n\cdot\Xi\cdot ?\quad}>>
\bK(\sB,n)\thf\,.
\endCD$$
\msk
\example{4.1.1. Summary} $\bL\bA^h\bul(X,n)$  is a
spectrum  defined for any space 
$X$ over $B\!\G$ and any $n\ge0$. It is the 
homotopy pullback of a diagram 
\msk
$$\CD
 @.  \bA(X,n)\hf \\
@.        @VV\cap V \\
\bL^h\bul(X)@>{\quad t^n\cdot\Xi\cdot ?\quad}>>
\bA(X,n)\thf
\endCD$$
\msk
in which $\bL^h\bul(X)$ denotes the 0--connected $L$--theory
spectrum of $\ZZ\pi_1(X)$ with decoration $h$~, and 
$\bA(X,n)=\SS^n_!\wedge\bA(X)$.  Hence there are
homotopy fiber sequences
$$\aligned
\bA(X,n)\ho@>>>&\bL\bA^h\bul(X,n)@>>>\bL\bul^h(X)~, \\
\bL\bA^h\bul(X,n-1)@>>>&\bL\bA^h\bul(X,n)
@>v>>\bA(X,n)\,.\endaligned\tag4.1.2$$
\endexample 
\bsk
\head 4.2. $LA$--theory and $h$--structure Spaces \endhead
In the following theorem, we mean by $(\bL\bA^h\bul)\upr(\text{---},n)$
and $(\bL\bA^h\bul)\lpr(\text{---},n)$
domain and homotopy fiber, respectively, 
of the assembly transformation for the homotopy 
functor $\bL\bA^h\bul(\text{---},n)$ on $\sW_{B\!\G}$. The manifold $M$ 
becomes an object in  $\sW_{B\!\G}$ by means of the classifying map for
$\nu(M)$. 
\msk
\proclaim{4.2.1. Theorem} For compact $M^n$
there exists a homotopy commutative square 
with highly connected (see {\rm 4.2.2}) vertical arrows
$$\CD
\sS(M)@>\nabla>>\sS(\tau) \\
@VVV   @VVV  \\
\Omega^{\infty+n}(\bL\bA^h\bul)\lpr(M,n)@>\text{forget}>>
\Omega^{\infty+n}(\bL\bA^h\bul)\upr(M,n)\,.
\endCD$$ 
\endproclaim
\msk
\definition{4.2.2. Details} The right-hand vertical arrow in 4.2.1 is 
$(j+2)$--connected if $j$ is in the
smooth concordance stable range for $\DD^n$ and $j\le n-2$.
The left--hand vertical arrow in 4.2.1 induces a bijection on $\pi_0$. Each 
component of $\sS(M)$ determines a homeomorphism class of manifolds $N$ 
homotopy equivalent to $M$. If $j$ is in the topological concordance stable 
range for $N$, then the left--hand vertical arrow in 4.2.1 restricted to that 
component (and the corresponding component of the codomain) is
$(j+1)$--connected.
\enddefinition
\msk
\example{4.2.3. Remark} There is an $s$--decorated version of 4.2.1, in
which  the Space of $s$--structures $\sS^s(M)$ replaces $\sS(M)$ and
$\bL\bA^s\bul$ replaces $\bL\bA^h\bul$. To define $\bL\bA^s\bul$
use $L$--theory and algebraic $K$--theory of spaces with an
$s$--decoration in
\S4.1.  The homotopy groups of $\bL\bA^s\bul(X,n)$ differ from
those of
$\bL\bA^h\bul(X,n)$  only in dimensions $\le n$. 
\msk
A ${\langle-\infty\rangle}$--decorated
version of 4.2.1 exists, but does not give anything new since 
the homotopy groups of 
$$\bL\bA^{\langle-\infty\rangle}\bul(\text{---},n)$$
differ from those of $\bL\bA^h\bul(\text{---},n)$ only  in dimensions $<n$.
Nevertheless, it is good to have this in mind when making 
the comparison with 2.5.1 (next remark).
\endexample
\msk
\example{4.2.4. Remark} Theorems 4.2.1 and 2.5.1 are compatible: the
commutative squares in 4.2.1 and 2.5.1  are opposite faces of a
commutative cube. Of the 12 arrows in  the cube, the four not mentioned
in 4.2.1 or 2.5.1 are inclusion maps (top) and forgetful maps (bottom). 
Also, the $s$--version of 4.2.1 is compatible with 2.3.1 in the same sense.
\endexample
\msk
\example{4.2.5. Remark} The left--hand column of the 
diagram in 4.2.1 matches the left--hand column of the diagram
in 3.2.1.  In detail:  there exists a commutative diagram 
with highly connected vertical arrows
\msk
$$\CD 
\sH^{\infty}(M)@<w<<\sS(M) \\
@VVV  @VVV \\
\Omega\8(\bA\lpr(M))@<v\lpr<<\Omega^{\infty+n}(\bL\bA^h\bul)\lpr(M,n)
\endCD$$
\msk
where $w$ is the Whitehead torsion map of 1.5 and 
$v\lpr$ is induced by $v$ of \thetag{4.1.2}. 
\endexample
\bsk
\head 4.3.  Special features of the smooth case  \endhead 
In this section we assume that $M^n$ is smooth. We use 
the notation of 4.2.1 and 2.5.1. By 4.2.4,
there is a commutative diagram 
\msk
$$\CD\sS(\tau)@>\subset>> \bigcup_k\sS(\tau\oplus\ep^k) \\
@VVV     @VVV \\
\Omega^{\infty+n}(\bL\bA^h\bul)\upr(M,n)
@>>>\Omega^{\infty+n}(\bL^h\bul)\upr(M,n)\,.\endCD\tag4.3.1$$
\msk
We can write the resulting map of horizontal homotopy 
fibers in the form 
$$\CD
u\sS(\tau) \\
@VVV \\
\Omega^{\infty+n}(\bA\upr(M,n)\ho)
\endCD\tag4.3.2$$
where the prefix $u$ indicates {\it unstable} structures. It is
highly connected~; details as in 4.2.2. Our goal here is to 
describe a smooth analog of \thetag{4.3.2}. 
\msk
\proclaim{4.3.3. Proposition} For a space $X$ over $B\!\Or$~, 
the map $\Sigma\8(X_+)\to \bA(X)$ 
of {\rm 3.2.2} has a canonical refinement to a $\ZZ/2$--map. 
\endproclaim
\msk
\demo{Remarks, Notation} As in \S4.1, we work in ``naive''  
stable $\ZZ/2$--homotopy theory.  That is, whenever 
we see a $\ZZ/2$--map $\bY'\to\bY$ which is an
ordinary homotopy equivalence, we are allowed 
to replace $\bY$ by $\bY'$. 
 
The involution on $\bA(X)$ needed here is
determined by $X\to B\!\Or\hookrightarrow B\!\G$ as in \S4.1.
The involution on $\Sigma\8(X_+)$ that we have in mind is 
as follows. For simplicity, assume that $X$ is a compact
CW--space~; then the reference map $X\to B\!\Or$ factors
through 
$B\!\Or(k)$ for some $k$. Let $\gamma^k$ be the vector bundle 
on $X$ pulled back from $B\!\Or(k)$, and let $\eta^{\ell}$ be a 
complementary vector bundle on $X$~, so
that $\gamma\oplus\eta$ is  trivialized. Now, to see the 
involution, replace $\Sigma\8(X_+)$ by the homotopy 
equivalent 
$$\Omega^{\ell}_!\Omega^k
\Sigma\8\Sigma^{\gamma}\Sigma^{\eta}_!(X_+)$$
where $\Sigma^{\eta}(X_+)$ and
$\Sigma^{\gamma}\Sigma^{\eta}(X_+)$ are the Thom spaces of
$\eta$ and $\gamma\oplus\eta$~, respectively.
Subscripts $!$ indicate that $\ZZ/2$ acts on loop or suspension
coordinates by scalar multiplication with $-1$.  Compare \S4.1. 
--- We abbreviate
$$\Sigma\8(X_+~,n):=\SS^n_!\wedge\Sigma\8(X_+)\,.$$
\enddemo
\msk
While the map $\Sigma\8(X_+)\to \bA(X)$ of 3.2.2 can be 
described in algebraic $K$--theory terms, 
including the algebraic $K$--theory of finite sets over $X$~, our
proof  of 4.3.3 is not entirely $K$--theoretic. It uses 3.2.1 to 
interpret the involution on $A(X)$ geometrically. 
\msk
\proclaim{4.3.4. Theorem} There is a commutative diagram
with highly connected vertical arrows
(and lower row resulting from {\rm 4.3.3})
\msk
$$\CD
u\sS(\tau) @<<<u\sS_d(\tau)\\
@VVV  @VVV \\
\Omega^{\infty+n}(\bA\upr(M,n)\ho)
@<<< \Omega^{\infty+n}((\Sigma\8(X_+~,n))\ho)\,.
\endCD$$
\endproclaim
\msk
\demo{Remarks} The right--hand column of this diagram 
is $(n-2)$--connected. It is essentially an old construction 
due to Toda and James; see \cite{Jm} for references.  

Something should be said about compatibility between 4.3.4
and 4.2.5, but we will leave it unsaid. 
\enddemo
\msk
\definition{4.3.5. Remark} In calculations involving 4.3.4 the 
concept of {\it stabilization} is often useful. Stabilization
is a way to make new homotopy functors 
on $\sW_Y$ (spaces over $Y$) from old ones.
Idea: Given a homotopy functor 
$\bJ$ from $\sW_Y$ to spectra, and $X$ in $\sW_Y$, let $s\bJ(X)$ 
be the homotopy colimit of the $\Omega^n(\bJ(X\times\SS^n)/\bJ(X))$ for
$n\ge0$.  There is a natural transformation $\bJ(X)\to s\bJ(X)$
induced by $x\mapsto (x,-1)\in X\times\SS^0$ for $x\in X$. 
The formalities are much as in \S2.4, even though the  result is
quite different. The main examples for us are these:
\roster
\itb Take $\bJ(X)=\bA(X)$ for $X$ in $\sW_{B\!\Or}$. Then 
$s\bA(X)\simeq\Sigma\8(\Lambda X_+)$ where $\Lambda X$ 
is the free loop space. (See \cite{Go1} for details.) Hence
$(s\bA)\upr(X)\simeq\Sigma\8(X_+)$.
\itb Take $\bJ(X)=\bL\bul^h(X)$ for $X$ in $\sW_{B\!\Or}$. 
Then $s\bL\bul^h(X)$ is contractible  for all $X$ by the
$\pi$--$\pi$--theorem. Hence $(s\bL\bul^h)\upr(X)$ 
is also contractible. 
\endroster
One can use these facts to split the lower row in 4.3.4,
up to homotopy. See 6.5 for another application.  
\enddefinition
\bsk
\head  {\bf 5.  Geometric structures on fibrations} \endhead
\head 5.1. Block bundle structures \endhead 
Here we address the following question. Given a fibration $p$ on a Space $B$ 
whose fibers are Poincar\'e duality spaces of formal dimension $n$, can we 
find a block bundle $p_0$ on $B$ with closed manifold fibers, fiber
homotopy  equivalent to $p$~? For earlier work on this problem,
see \cite{Qun4}, \cite{Qun1}. We combine this with ideas from
\cite{Ra4} and \cite{Ra2}.  
\msk
The {\it block $s$--structure Space} $\wdt{\sS}^s(X)$ 
of a simple Poincar\'e duality space $X$ of formal dimension $n$ 
(alias {\it finite} Poincar\'e space, \cite{Wa1, \S2}~)  is defined literally
as in the case of a closed manifold. Any  simple homotopy equivalence
$M^n\to X$, where
$M^n$ is a closed  manifold, induces a homotopy equivalence 
$$\wdt{\sS}^s(M)@>>>\wdt{\sS}^s(X)$$
and $\wdt{\sS}^s(M)$ was described in $L$--theoretic
terms in 2.3.1.  But such a
homotopy equivalence $M\to X$ might not exist, and even if it does, we
might want to see  an $L$--theoretic description of the block 
$s$--structure Space of $X$ which does not use a choice 
of homotopy equivalence $M\to X$. 
\msk
Ranicki \cite{Ra4}, \cite{Ra2, \S17} associates to a simple Poincar\'e
duality space 
$X$ of formal dimension $n$ its {\it total surgery obstruction}, a point 
$$\partial\sigma^*(X)\in \Omega^{\infty+n-1}((\bL\bul^s)\lpr(X))\,.$$
The element has certain naturality 
properties. For example, a homotopy equivalence $g:X\to Y$ 
determines a path in $\Omega^{\infty+n-1}((\bL\bul^s)\lpr(X))$
from $g_*\partial\sigma^*(X)$ to 
$\partial\sigma^*(Y)$;
 more  later. 
\msk
\proclaim{5.1.1. Theorem} $\wdt{\sS}^s(X)$ is (naturally) 
homotopy equivalent to the space of paths from 
$\partial\sigma^*(X)$ to the base point in 
$\Omega^{\infty+n-1}((\bL\bul^s)\lpr(X))$.
\endproclaim
\msk
Note that any choice of (base) point in the space of paths in 5.1.1
leads to an identification of it with
$\Omega^{\infty+n}((\bL\bul^s)\lpr(X))$, up to homotopy
equivalence. In this way, we recover the result 

$$\wdt{\sS}^s(M)\simeq \Omega^{\infty+n}((\bL\bul^s)\lpr(M))\,.$$
The advantages of 5.1.1 become clearer when it is applied
to {\it families}, that is, fibrations $p:E(p)\to B$ whose fibers are 
Poincar\'e  spaces of formal dimension $n$.  
(One must pay attention to simple 
homotopy types, so we assume that $B$ is connected 
and $p$ is classified by a map $B\to B\!\G^s(X)$ for some simple Poincar\'e duality 
space $X$ as above.) Given such a fibration $p:E(p)\to B$, we obtain 
an associated fibration $q:E(q)\to B$ with fiber 
$$\Omega^{\infty+n-1}((\bL\bul^s)\lpr(p^{-1}(b)))$$
over $b\in B$, 
and a section $\partial\sigma^*(p)$ of $q$ selecting the
total surgery obstruction $\partial\sigma^*(p^{-1}(b))$ in
$q^{-1}(b)$, for $b\in B$. The fibers of $q$ are 
infinite loop spaces, so we also have a zero section. (Technical point: For these 
constructions it is convenient to assume that $B$ is a simplicial complex, and to 
apply a suitable $(n+2)$--ad version of 5.1.1 to $E(p)|\sigma$ for each $n$--simplex $\sigma$ 
in $B$.) We say that $p$ {\it admits a block bundle structure} if 
the classifying map $B\to B\!\G^s(X)$ lifts to a map
$$B\to B\wdt{\TOP}(M)$$ for some closed manifold  $M$ equipped with a simple 
homotopy equivalence to $X$. 
\msk
\proclaim{5.1.2. Corollary} The fibration $p:E(p)\to B$ with 
Poincar\'e duality space fibers admits a block bundle
structure  if and only if $\partial\sigma^*(p)$ is vertically
nullhomotopic. 
\endproclaim
\msk
In the case $B=B\!\G^s(X)$ we can add the following:
\msk
\proclaim{5.1.3. Corollary} Let $p:E(p)\to B\!\G^s(X)$ be the 
canonical fibration with fibers $\simeq X$. There is a 
cartesian square
$$\CD
\coprod\limits_M B\wdt{\TOP}(M)@>>> B\!\G^s(X) \\
@VVV   @VV\text{total surgery obstruction section}V  \\
B\!\G^s(X)@>\text{zero section}>> E(q)
\endCD\qquad$$
where $M$ runs through a maximal set of pairwise 
non--homeomorphic closed $n$--manifolds in the simple 
homotopy type of $X$. 
\endproclaim
\msk
A few words on how $\partial\sigma^*(X)$
is constructed:  Ranicki creates a homotopy functor 
$\bV\bL_s\ubul$ on spaces over $\RR P\8$, and a natural transformation 
$\bL^s\bul\to \bV\bL_s\ubul$ with the property that 
$$\CD
(\bL^s\bul)\upr(X) @>>> (\bV\bL_s\ubul)\upr(X) \\
@VV\text{assembly} V    @VV\text{assembly} V \\
\bL^s\bul(X)@>>> \bV\bL_s\ubul(X)
\endCD$$
is cartesian for any $X$ over $\RR P\8$.  (The functor $\bV\bL_s\ubul$
we have in mind is defined in \cite{Ra, \S15},  and we should 
really call it $\bV\bL_s\ubul\langle1/2\rangle$ to conform with 
Ranicki's notation.) Any finite
Poincar\'e duality  space $X$ of formal dimension $n$ determines an
element 
$\sigma^*(X)\in \Omega^{\infty+n}(\bV\bL_s\ubul(X))$, the 
{\it visible symmetric signature} of $X$. The image of $\sigma^*(X)$ 
under the boundary map 
$$\Omega^{\infty+n}(\bV\bL_s\ubul(X))
\la \Omega^{\infty+n-1}((\bV\bL_s\ubul)\lpr(X)) 
\quad\simeq\quad \Omega^{\infty+n-1}((\bL^s\bul)\lpr(X))$$
is the total surgery obstruction $\partial\sigma^*(X)$.  Therefore: 
$X$ is simple homotopy equivalent to a closed $n$--manifold if and only if 
the component of $\sigma^*(X)$ is in the image of the assembly 
homomorphism, 
$$\pi_n(\bV\bL_s\ubul)\upr(X)\to \pi_n\bV\bL_s\ubul(X)\,.$$
The functor $\bV\bL_s\ubul$
has a {\it ring structure}, that is, for $X_1$ and $X_2$ over $\RR P\8$ 
there is a multiplication 
$$\mu:\bV\bL_s\ubul(X_1)\wedge \bV\bL_s\ubul(X_2) \la
\bV\bL_s\ubul(X_1\times X_2),$$
with a unit in $\bV\bL_s\ubul(*)$. The visible symmetric 
signature is multiplicative:
$$\mu(\sigma^*(X_1),\sigma^*(X_2)) = \sigma^*(X_1\times X_2),$$
up to a canonical path, for Poincar\'e duality spaces $X_1$ and $X_2$. 
This property makes $\sigma^*$ useful (more useful than 
$\partial\sigma^*$ alone) in dealing with products, say, in 
giving a description along the lines of 5.1.1 of the product map 
$$\wdt{\sS}^s(X_1)\times\wdt{\sS}^s(X_2)@>>>
\wdt{\sS}^s(X_1\times X_2)\,.$$
\bsk
\head 5.2. Fiber bundle structures \endhead
Here the guiding question is: Given a fibration $p$ over some space $B$, with 
fibers homotopy equivalent to finitely dominated CW--spaces, does there exist 
a bundle $p_0$ on $B$ with compact manifolds as fibers, fiber
homotopy equivalent to $p$~? We do not assume that the fibers of
$p$ satisfy  Poincar\'e duality. We do not ask that 
the fibers of $p_0$ be closed and we do not care what dimension 
they have.  See \cite{DWW}, \cite{DWWc} for all details. 
\msk 
Let $Z$ be a compact CW--space, equipped with 
a euclidean bundle $\xi$.  Let $\sT_n(Z,\xi)$ be
the  Space of pairs $(M,f,j)$ where $M^n$ is a compact 
manifold with boundary, $f:M\to Z$ is a homotopy equivalence,
and $j$  is a stable isomorphism $f^*\xi\to\tau(M)$. Let $\sT(Z,\xi)$
be the colimit of the $\sT_n(Z,\xi)$ under stabilization
(product with $I$). 
\msk
Any choice of vertex $(M,f,j)$ in $\sT(Z,\xi)$ leads to a homotopy 
equivalence from $\sT(M,\tau)$ to $\sT(Z,\xi)$. There is a homotopy  
fiber sequence
$$\bigcup_k\TOP(M\times I^k,\partial(M\times I^k))
\la \bigcup_k\G\tan(M_0\times \RR^k)
\la \sT(M,\tau)$$
where $M_0=M\minus\partial M$, and 
$\G\tan(\dots)$ refers to homotopy automorphisms $f$
of $M_0\times \RR^k$ covered by isomorphisms 
$\tau(M_0\times\RR^k)\to f^*(\tau(M_0\times\RR^k))$.
By 1.2.1 we may write  $\bigcup_k\TOP(M_0\times \RR^k)$
instead of $\bigcup_k\G\tan(M_0\times \RR^k)$. 
Now 1.2.3 implies 
$$\Omega \sT(M,\tau)\,\simeq\,\Omega\sH^{\infty}(M)$$
and suggests $\sT(M,\tau)\,\simeq\,\sH^{\infty}(M)$. This is 
easily confirmed with the methods of \S1.2. Summarizing
these observations: any choice of vertex $(M,f,j)$ in $\sT(Z,\xi)$ leads to 
a homotopy equivalence $\sT(Z,\xi)\to\sH^{\infty}(M)$. 
Furthermore, $\sH^{\infty}(M)\simeq \Omega\8\bA\lpr(M)
\simeq\Omega\8\bA\lpr(Z)$ by 3.2.1.
\msk
Again, we might want to see a description of $\sT(Z,\xi)$ 
in terms of $\bA(Z)$ which does not depend on a choice of base point 
in $\sT(Z,\xi)$. To get such a description we proceed very much as 
in \S5.1, by associating to $Z$ a {\it characteristic}
$$\chi(Z)\in\Omega\8\bA(Z),$$
analogous to Ranicki's 
$\sigma^*(X)\in\Omega^{\infty+n}\bV\bL\ubul_s(X)$ for a
simple Poincar\'e duality space $X$ of formal dimension $n$. The
element 
$\chi(Z)$ is the image of the object/vertex 
$$\phantom{xxxxxxxxxxxxxxxxxxxxx}\SS^0\times Z
 \overset r\to{\overset{\dsize\la}\to{\underset i\to\longleftarrow}}
Z\qquad\quad (r(x,z)=z~,\,\, i(z)=(1,z))$$
 in $\sR(Z)$ under the inclusion 
$|w\sR(Z)|\hookrightarrow\Omega\8\bA(Z)$ mentioned in \S3.1. 
If $Z$ is connected, then  
the component of $\chi(Z)$ in $\pi_0\Omega\8\bA(Z)\cong\ZZ$ is 
the Euler characteristic of $Z$.
\msk
\proclaim{5.2.1. Theorem} $\sT(Z,\xi)$ is (naturally) homotopy 
equivalent to the homotopy fiber of the assembly map 
$\Omega\8(\bA\upr(Z))\la \Omega\8\bA(Z)$ over the point $\chi(Z)$.
\endproclaim
\msk
Note that the $A$--theoretic expression for $\sT(Z,\xi)$ does not 
depend on $\xi$. 
Again, {\it naturality} in 5.2.1 is a license to apply the 
statement to families.  Let $p:E\to B$ be a fibration where the
fibers $E_b$ are homotopy  equivalent to compact CW--spaces. 
Let $\Omega\8\bA(p)$ and 
$\Omega\8\bA\upr(p)$ be the associated fibrations on $B$ with 
fibers $\Omega\8\bA(E_b)$ and $\Omega\8\bA\upr(E_b)$,
respectively,  over a point $b\in B$. The rule $b\mapsto\chi(E_b)$ 
defines a section of $\Omega\8\bA\upr(p)$ which we call $\chi(p)$.
See \cite{DWW} for explanations regarding  the {\it continuity} of this construction.
Assembly  gives a map over $B$ from the total space of
$\Omega\8\bA\upr(p)$ to that of
$\Omega\8\bA(p)$. 
\msk
\proclaim{5.2.2. Corollary} The fibration $p:E\to B$ is fiber 
homotopy equivalent to a bundle with compact manifolds as fibers 
if and only if the section $\chi(p)$ of $\Omega\8\bA(p)$ lifts (after a
vertical  homotopy) to a section of $\Omega\8\bA\upr(p)$.
\endproclaim
\msk
We leave it to the reader to state an analog of 5.1.3, and turn instead
to the smooth case. Suppose that $\xi$ is a {\it vector bundle}
over $Z$. There is then a smooth variant $\sT_d(Z,\xi)$ of 
$\sT(Z,\xi)$.  Any choice of base vertex $(M,f,\xi)$ 
in $\sT_d(Z,\xi)$  leads to a homotopy equivalence
$$\sT(Z,\xi)\simeq \sH^{\infty}_d(M)\,.$$
Remember now (3.2.2) that $\sH^{\infty}_d(M)$ can be $A$--theoretically 
described as the homotopy fiber of some map 
$\eta:\Omega\8\Sigma\8(M_+)\to \Omega\8\bA(M)$. The map 
$\eta$ is, in homotopy invariant terms,  $\Omega\8$ of the composition 
$$M_+\wedge\bS^0@>\id\wedge u>>M_+\wedge\bA(*)
\quad\simeq\quad\bA\upr(M)\la\bA(M)$$ 
where $u:\bS^0\to\bA(*)$ is the unit of the ring spectrum $\bA(*)$. 
\msk
\proclaim{5.2.3. Theorem} $\sT_d(Z,\xi)$ is (naturally) homotopy 
equivalent to the homotopy fiber of 
$\eta:\Omega\8\Sigma\8(Z_+)\la \Omega\8\bA(Z)$ over the point 
$\chi(Z)$.
\endproclaim
\msk
Returning to the notation and hypotheses of 5.2.2, we are compelled 
to introduce yet another fibration $\Omega\8\Sigma\8_+p$ on $B$, with
fiber $\Omega\8\Sigma\8(E_b)_+$ over $b\in B$.
\msk
\proclaim{5.2.4. Corollary} The fibration $p:E(p)\to B$ is fiber 
homotopy equivalent to a bundle with smooth compact manifolds as fibers 
if and only if the section $\chi(p)$ of $\Omega\8\bA(p)$ lifts (after a
vertical  homotopy) to a section of $\Omega\8\Sigma\8_+p$.
\endproclaim
\msk
\remark{Remarks} In 5.2.4, {\it bundle with smooth compact 
manifolds as fibers} means, say in the case where $B$ is connected, 
a fiber bundle with fibers $\cong M$ where $M$ is smooth compact, 
and structure group $\DIFF(M,\partial M)$. 

Corollary 5.2.4 is closely related to something we shall discuss 
in \S6.7: the Riemann--Roch theorem of 
\cite{BiLo}, see also \cite{DWW}. 
\endremark
\msk
\proclaim{5.2.5. Corollary} Let $p:E\to B$ be a fibration with 
fibers homotopy equivalent to compact CW--spaces. If $Y$ 
is any compact connected CW--space of Euler characteristic $0$,
then  the composition $pq:Y\times E\to B$ (where $q:Y\times E\to E$
is the projection) is fiber homotopy equivalent to a bundle 
with smooth compact manifolds as fibers. 
\endproclaim
\msk
\demo{Proof} A suitable product formula implies that 
$\chi(pq)$ is vertically homotopic to $\chi(Y)\times\chi(p)$.
We saw earlier that the component of $\chi(Y)$ in 
$\pi_0\bA(Y)\cong\ZZ$ is the Euler characteristic of $Y$. \qed
\enddemo

Statements 5.2.1--5 can easily be generalized 
to the case where $Z$ is a {\it finitely dominated} CW--space. But
it is then necessary to use a variant $\bA^p(Z)$ of $\bA(Z)$ with a
larger $\pi_0$ isomorphic to $K_0(\ZZ\pi_1(X))$. Then 
$\chi(Z)$ in $\Omega\8\bA^p(Z)$ is defined, and 5.2.1 and 5.2.3
remain  correct as stated.  In this more general formulation, 
5.2.1 includes Wall's theory \cite{Wa2} of the finiteness
obstruction.

Casson and Gottlieb \cite{CaGo} established 5.2.5 in the case 
$Y=(\SS^1)^n$~, with a large $n$ depending on $p:E\to B$. 
\bsk
\head  {\bf 6. Examples and Calculations} \endhead
\example{6.1.  Smooth automorphisms of disks} 
The smooth version of 
\thetag{1.2.2} gives a  homotopy fiber sequence
$$\DIFF(\DD^n,\SS^{n-1})\la\DIFF(\RR^n)\to \sH_d(\SS^{n-1})$$
where $\sH_d$ is a Space of differentiable $h$--cobordisms.
The composition of group homomorphisms
$\Or(n)\hookrightarrow \DIFF(\DD^n,\SS^{n-1})\to \DIFF(\RR^n)$
is a homotopy equivalence, so that 
$$\DIFF(\DD^n,\SS^{n-1})\,\,\simeq\,\,
\Or(n)\times \Omega\sH_d(\SS^{n-1})\,.$$
By 1.3.4 and 3.2.2, there is a
map from $\Omega\sH_d(\SS^{n-1})$ to 
$\Omega^{\infty+2}\Whd(\SS^{n-1})$ 
which is an isomorphism on $\pi_j$ for
$j<\phi(n-1)$, where $\phi(n)$ is the
minimum of $(n-4)/3$ and $(n-7)/2$.  Further,  the map
$\Whd(\SS^{n-1})\to \Whd(*)$ induced by 
$\SS^{n-1}\to *$ is
approximately 
$2n$--connected \cite{Wah2} and the rational homotopy groups of 
$\Whd(*)$ in dimensions $>1$ are those of $\bK(\QQ)$, which 
are known \cite{Bo}.  Therefore:
$$\pi_j\DIFF(\DD^n,\SS^{n-1})\,\otimes\,\QQ\cong\left\{\aligned 
&\QQ\oplus\QQ\qquad\text{if $0< j<\phi(n-1)$ and  $4\mid j+1$}
\\ &0 \qquad\qquad\text{if $0< j<\phi(n-1)$ and  $4\nmid j+1$.}
\endaligned\right.$$
For a calculation of the rational homotopy groups of 
$\DIFF(\DD^n)$ in the concordance stable range, following 
Farrell and Hsiang \cite{FaHs}, we note $\DIFF(\DD^n)
\simeq\Omega\sS_d(\DD^n)$ and use the smooth version of 
1.5.2, which gives (at odd primes and in the concordance stable 
range)
$$\sS_d(\DD^n)\simeq\wdt{\sS}_d(\DD^n)\times
\wdt{\DIFF}(\DD^n)/\DIFF(\DD^n)\simeq
\wdt{\sS}_d(\DD^n)\times \Omega\8(\bH_d(\DD^n)\ho)\,.$$
Here $\wdt{\sS}_d(\DD^n)\simeq\Omega^n(\TOP/\Or)$, which 
is rationally trivial, and $\bH_d(\DD^n)$ is homotopy 
equivalent to $\Omega\Whd(*)$, so
3.2.2 and \cite{Bo} give
$$\pi_j\bH_d(\DD^n)\otimes\QQ\cong \left\{\aligned &\QQ
\qquad\text{if $4|j$} \\ &0\qquad\text{otherwise}
\endaligned\right. $$
provided $0<j\le\phi(n)$. The canonical involution on $\bH_d(\DD^n)$ acts
trivially  on these rationalized homotopy groups if $n$ is odd, and 
nontrivially if $n$ is even. Therefore, if $0\le j<\phi(n)$, then 
\cite{FaHs}
$$\pi_j\DIFF(\DD^n)\,\otimes\,\QQ\cong\left\{\aligned 
&\QQ\qquad\text{$n$ odd and  $4\mid j+1$}
\\ &0 \qquad\qquad\text{otherwise.}
\endaligned\right.$$
Beware that Farrell and Hsiang write $\DIFF(\DD^n,\partial)$ for 
our $\DIFF(\DD^n)$.
\endexample
\msk
\example{6.2. Smooth automorphisms of spherical space forms} Let $M^n$ be smooth 
closed orientable, with universal cover $\simeq\SS^n$, where $n\ge5$.
Hsiang and Jahren \cite{HsiJ} calculate $\pi_*\DIFF(M)\otimes\QQ$ in the smooth 
concordance stable range, assuming that $n$ is odd. They begin with the 
observation that $\pi_j\G^s(M)$ is finite for all $j$. Therefore 
$$\pi_j\DIFF(M)\otimes\QQ \cong \pi_{j+1}\sS^s_d(M)\otimes\QQ$$ 
for $j>0$. By the smooth versions of 1.5.2 and 1.4.4
we have a splitting 
$$\sS^s_d(M)\,\simeq\,\wdt{\sS}^s_d(M)\times\Omega\8(\bH^s_d(M)\ho)$$
at odd primes, in the concordance stable range. Therefore, rationally, 
$$\pi_j\DIFF(M)\cong \pi_{j+1}\sS^s_d(M)\cong
\pi_{j+n+2}\wdt{\bL}^h\bul(M)\oplus\pi_{j+1}\bL^h\bul(*)\oplus\pi^-_{j+1}\bH^s_d(M)
$$
for $0<j<\phi(n)$, where  $\wdt{\bL}^h\bul$ 
is the reduced $L$--theory and $\pi^-_{j+1}\bH^s_d(M)$ is the 
quotient of $\pi_{j+1}\bH^s_d(M)$ by the fixed subgroup of the 
$\ZZ/2$--action. This is the Hsiang--Jahren
result.  Rationally, the multisignature homomorphisms on 
$\pi_*\bL^h\bul(M)$ are   isomorphisms \cite{Wa3}. Rationally,
$\pi_*(\bH^s_d(M))\cong\pi_{*+1}\bK(\QQ\pi_1(M))$ for $0<*<n-1$.  The
calculation of $\pi_*\bK(\QQ\pi)\otimes\QQ$ for a finite group $\pi$  can 
often be accomplished with \cite{Bo}, certainly  in the case where
$\pi_1(M)$ is commutative. 
\endexample
\msk
\example{6.3.  Automorphisms 
of negatively curved
manifolds} Let $M^n$ be smooth, closed,  connected, 
with a Riemannian metric of  sectional curvature $<0$. In the course
of their proof of the Borel conjecture for such $M$~, Farrell and
Jones \cite{FaJo1}, \cite{FaJo5} show that 
$\alpha:(\bL\bul)\upr(M)\to\bL\bul(M)$ is a homotopy 
equivalence (with a $4$--periodic definition of $\bL\bul(M)$,
decoration $s$ or $h$). With our 0--connected definition 
of $\bL\bul(M)$, it is still true that 
$$\gather \Omega^{\infty+n}((\bL\bul)\lpr(M))\,\simeq\,* \\
\Omega^{\infty+n}((\bL\bA^h\bul)\lpr(M,n))\,\simeq\,
\Omega^{\infty+n}(\bA\lpr(M,n)\ho)\,.\endgather$$
For a simple closed geodesic $T$ in 
$M$, let $T\sha$ be the ``desingularization''  of $T$, 
so that $T\sha\cong\SS^1$. Farrell and Jones also show that the map 
$$\bigvee_T \bA\lpr(T\sha)\to \bA\lpr(M)$$
induced by $T\sha\to T\hookrightarrow M$ for all simple 
closed geodesics $T$ in $M$ is a homotopy equivalence \cite{FaJo3}; see also 
\cite{FaJo4}, \cite{FaJo2} for extensions. 
Now  there is a {\it fundamental theorem} in the algebraic
$K$--theory of spaces:  The assembly from 
$\SS^1\wedge\bA(*)\simeq\bA\upr(\SS^1)$  to $\bA(\SS^1)$
is a split monomorphism in the homotopy category
\cite{KVWW1}, \cite{KVWW2} and its mapping cone splits up to 
homotopy into two copies 
of a spectrum $\Nil_{\bA}(*)$. Under any of the involutions 
constructed by the method of \S4.1, these two copies 
are interchanged. Therefore 
$\Omega^{\infty+n}((\bL\bA\bul)\lpr(M,n))\simeq
\Omega^{\infty+1}(\bigvee_T\Nil_{\bA}(*))$
and we get from 4.2.1 a map 
$$\sS(M)\to \Omega^{\infty+1}(\bigvee_T\Nil_{\bA}(*))$$
which is approximately $(n/3)$--connected (see 4.2.2). 
It is known that $\Nil_{\bA}(*)$ is rationally  trivial and
$1$--connected, but $\pi_2\Nil_{\bA}(*)\ne 0$. 
See \cite{HaWa},\cite{Wah1}.
>From \cite{Dun}, \cite{BHM}, one has a homological algebra description
of $\Nil_{\bA}(*)$, as explained in \cite{Ma1, 4.5} and \cite{Ma2, \S5}.
But the homological algebra is over the ring
spectrum $\bS^0$ and it is not considered easy.---
>From the fiber sequence $\TOP(M)\to \G(M)@>>>\sS(M)$
we get
$\pi_j\TOP(M)\cong\bigoplus_T\pi_{j+2}\Nil_{\bA}(*)$
if $1<j< \phi(n)$~, and an exact sequence 
$$\align
\bigoplus_T\pi_3\Nil_{\bA}(*)\rightarrowtail&\pi_1\TOP(M)\to
\text{center}(\pi_1(M)) \\
\la &\bigoplus_T\pi_2\Nil_{\bA}(*)\to\pi_0\TOP(M)
\twoheadrightarrow\text{Out}(\pi_1(M)).
\endalign$$
\endexample 
\msk
\example{6.4. The $h$--structure Space of $\SS^n$, for $n\ge5$} 
By 4.2.1 and 4.2.2 there is a commutative square 
$$\CD
\sS(\SS^n) @>>>\Omega^{\infty+n}((\bL\bA\bul)\lpr(\SS^n,n))  \\
@VV\cap V      @VVV \\
\wdt{\sS}(\SS^n)
@>\simeq>>\Omega^{\infty+n}((\bL^s\bul)\lpr(\SS^n)),
\endCD$$
where the  top horizontal arrow is highly connected.
Here $\Omega^n(\bL^s\bul)\lpr(\SS^n)$ simplifies to $\bL\bul(*)$,
and $\Omega^n(\bL\bA\bul)\lpr(\SS^n,n)$ has an analogous simplifying
map (not a homotopy  equivalence, but highly connected) to
$\bL\bA\bul(*,n)$. Summarizing: there is a homotopy commutative
square 
$$\CD 
\sS(\SS^n)@>>>\Omega\8\bL\bA\bul(*,n)\\
@VV\cap V @VV \text{forget}  V \\
\wdt{\sS}(\SS^n) @>\simeq >>\Omega\8\bL\bul(*) \endCD$$
and a homotopy fiber sequence 
$\bL\bA\bul(*,n)\to \bL\bul(*)@>\,\,\beta\,\,>>
\SS^1\wedge(\bA(*,n))\ho$ from \thetag{4.1.2}.  Calculations 
\cite{WWp} using connective $K$--theory $\bo$ as a substitute  for 
$\bA(*)$, via 
$\bA(*)\to\bK(\ZZ)\to\bo$,
show that $\beta$ detects all elements in $\pi_{n+q}\bL(*)$ whose signature is 
not divisible by $2a_q$ if $4$ divides $n+q$, and by 
$4a_q$ if $4$ divides both $n+q$ and $q$.  Here $a_q=1,2,4,4,8,8,8,8$
for $q=1,2,\dots,8$ and $a_{q+8}=16a_q$; the numbers $a_q$ are 
important in the theory of Clifford modules \cite{ABS}.
Consequently, if $4$ divides $n+q$, then the image of  the inclusion--induced  homomorphism 
$$\pi_{n+q}\sS(\SS^n)\la \pi_{n+q}\wdt{\sS}(\SS^n)\cong 8\ZZ$$
is contained in $2a_q\ZZ$ if $4$ does not divide $q$, and in $4a_q\ZZ$
if  4 does divide $q$. Note the similarity of this 
statement with \cite{At, 3.3}, \cite{LM, IV.2.7}, and \cite{Tho, Thm.14}.
\endexample 
\msk
\example{6.5. Obstructions to unblocking smooth block
automorphisms} \linebreak One of the main points of  \S4.2 and the
introduction  to \S4 is a homotopy commutative diagram, for
compact
$M^n$ with $n\ge5$,  
\msk
$$\CD
\Omega^{\infty+n+1}(\bL^s\bul)\lpr(M)
@>>>\Omega^{\infty+n}(\bA^s\lpr(M,n)\ho) \\
@VV\simeq V   @AAA \\
\Omega\wdt{\sS}^s(M)@>>> \wdt{\TOP}(M)/\TOP(M)
\endCD\tag6.5.1$$
\msk
in which the upper row is  the connecting map 
from the first of the two homotopy fiber sequences in \thetag{4.1.2}, 
$$\Omega\bL^s\bul(M)\to \bA^s(M,n)\ho~,\tag6.5.2$$
with ${}\lpr$ and decoration $s$ and $\Omega^{\infty+n}$ inflicted.
Modulo the identification 
$\bA^s\lpr(M)\simeq\bH^s(M)$ (the $s$--decorated version of 3.2.2),
the right--hand column  of \thetag{6.5.1} is the highly connected
map which we found at the end of \S1.4 using purely geometric
methods.  The lower row of \thetag{6.5.1} is  the connecting map
from the homotopy fiber sequence 
$$\wdt{\TOP}(M)/\TOP(M)\la \sS^s(M) \la \wdt{\sS}^s(M)\,.$$
One of the main points of \S4.3 is that much of \thetag{6.5.1} has a
smooth  analog, in the shape of a homotopy commutative diagram 
\msk
$$\CD
\Omega^{\infty+n+2}\bL^s\bul(M)
@>>>\Omega^{\infty+n+1}(\Whd^s(M,n)\ho) \\
@VVV   @AAA \\
\Omega\wdt{\sS}^s_d(M)@>>> \wdt{\DIFF}(M)/\DIFF(M)
\endCD\tag6.5.3$$
\msk
defined for smooth compact $M$. Here 
$\Whd^s(M)$ is the mapping cone of the map 
$\Sigma\8(M_+)\to \bA^s(M)$ discussed in 4.3.3 (except for a 
decoration $s$ which we add here), 
and $\Whd^s(M,n):=\SS^n_!\wedge\Whd^s(M)$.  The upper 
row in \thetag{6.5.3} is $\Omega^{\infty+n+1}$ of \thetag{6.5.2}
composed with the projection $\bA^s(M,n)\ho\la
\Whd^s(M,n)\ho$. Again, modulo  an identification
of $\Whd^s(M)$ with $\bH^s_d(M)$~, coming  from \thetag{3.2.3}, 
the right--hand column of
\thetag{6.5.3} is  a purely geometric construction going back to (the
smooth  version of) \S1.4.  It is highly connected. The
left--hand column of \thetag{6.5.3} is not a homotopy equivalence,
which makes the analogy a little imperfect.  We
arrive at \thetag{6.5.3} by first making full use of 
4.2.1 and 4.2.4 to produce a {\it framed}
version  of \thetag{6.5.1}, with upper left--hand vertex 
$\Omega^{\infty+n+2}\bL\bul^s(M)$~, upper right--hand vertex 
$\Omega^{\infty+n+1}(\bA(M,n)\ho)$~, 
and lower left--hand vertex equal to $\Omega$ of the homotopy
fiber of  
$$\nabla:\wdt{\sS}^s@>>>\wdt{\sS}(\tau)\,.$$

Now assume that the classifying map $M\to B\!\Or$ for the 
stable tangent bundle factors up to homotopy through an 
aspherical space. (For example,  this is the case if $M$ is stably
framed.) Stabilization arguments as in 4.3.5 show that 
then $\Whd^s(M,n)\ho$ splits off $\bA^s(M,n)\ho$~, up to 
homotopy. Also, the connecting map \thetag{6.5.2}
factors through the summand $\Whd^s(M,n)\ho$~, up to homotopy. 
It follows that elements in $L^s_k(\ZZ\pi)$ detected under 
\thetag{6.5.2} are also detected under 
$$L^s_k(\ZZ\pi)\la \pi_{k-n-2}(\wdt{\DIFF}(M)/\DIFF(M))~,\tag6.5.4$$
the homomorphism coming from \thetag{6.5.3}. To exhibit such 
elements in $L^s_k(\ZZ\pi)$~, we suppose in addition that
$M$ is orientable and $\pi=\pi_1(M)$ is finite. Then we have  the
multisignature homomorphisms
$$L_k^s(\ZZ\pi)\to L_k^p(\RR\pi)
\cong\bigoplus_V L_k^p(E_V)\tag6.5.5$$
where the direct sum is over a maximal set of pairwise
non--isomorphic irreducible  real representations $V$ of $\pi$, and
$E_V$ is the endomorphism ring of $V$, isomorphic to 
$\RR$, $\CC$ or $\HH$ equipped with a standard conjugation 
involution \cite{Le}, \cite{Wa3}. It is known that 
$L_k^p(E_V)\cong\ZZ$ if $4|k$, and also $L_k^p(E_V)\cong\ZZ$ 
if $2|k$ and $E_V\cong\CC$; otherwise $L_k^p(E_V)=0$. A calculation
similar  to that mentioned in 6.4 shows that an element in 
$L_k^s(\ZZ\pi)=\pi_k(\bL^s\bul(M))$ will be detected by the
homomorphism associated with \thetag{6.5.2}
if, for some irreducible real representation $V$ of $\pi$, the
$V$--component of its multisignature is  not divisible by 
$$\aligned &2a_{k-n}\qquad\quad\text{(assuming $4|k$ and
$E_V\cong\RR$)} \\ 
&2a^c_{k-n}\qquad\quad\text{(assuming $2|k$ and
$E_V\cong\CC$)}
\\ &a_{k-n+4}/4\qquad\text{(assuming $4|k$ and $E_V\cong\HH$)}
\endaligned$$
where $a^c_q=1$ if $q=1,2$ and $a^c_{q+2}=2a^c_q$ for $q>2$. 
If $4|n$ and $V$ is the trivial 1--dimensional representation, we can 
can do a little better:  the element will also be detected if the 
$V$--component of its multisignature is not divisible by 
$$4a_{k-n}\qquad\quad\text{(assuming $4|k$)}\,.$$
Now the multisignature homomorphisms \thetag{6.5.5} are of course
periodic in $k$ with period 4, and are  rational isomorphisms 
\cite{Wa1}. Therefore many elements 
in $L^s_*(\ZZ\pi)$ are indeed detected by \thetag{6.5.2}, 
and a fortiori by \thetag{6.5.4}. 
\endexample
\demo{Remark} This calculation can be viewed as a cousin of Rochlin's 
theorem. To make this clearer we switch from  block diffeomorphisms to
bounded diffeomorphisms, i\.e\. we look at 
$$
L^{\langle -\infty\rangle}_k(\ZZ\pi)\quad@>\qquad>>\quad
\pi_{k-n-2}(\DIFF^b(M\times\RR\8)/\DIFF(M))
\tag6.5.6$$
instead of \thetag{6.5.4}. Compare 2.5.1 and 4.2.3. The same calculations as
before show that  an element $x$  in the domain of \thetag{6.5.6} maps 
nontrivially if, for some irreducible real representation $V$ of $\pi$, 
the $V$--component of the multisignature of $x$ is not divisible by
certain  powers of 2, depending on $V$, $k$ and $n$, exactly as above. 
Note in passing that 
$$L^p_k(E_V)\cong L^{\langle -\infty\rangle}_k(E_V)$$
so that we can indeed speak of multisignatures as before.  
Specializing to $M=*$ and $4|k$, with $k>0$, we see that elements $x$ whose
signature  is not divisible by $4a_k$ are detected by \thetag{6.5.6}.  But in
the case $M=*$  we also have
$\DIFF^b(M\times\RR\8)\simeq\Omega(\TOP\!/\!\Or)$ and we may
identify \thetag{6.5.6} with the boundary map  in the long exact sequence
of homotopy groups associated with the  homotopy fiber sequence
$$\TOP\!/\!\Or\la \G\!/\!\Or\la \G\!/\!\TOP\,.$$
Therefore, if $4|k$ and $k>0$, the image of $\pi_k(\G\!/\!\Or)$ 
in $\pi_k(\G\!/\!\TOP)=8\ZZ$ is contained in $4a_k\cdot\ZZ$.  For $k=4$, 
this statement is (one form of) Rochlin's theorem.  For $k>4$, it is 
also well known as the 2--primary aspect 
of the Kervaire--Milnor work on homotopy 
spheres \cite{KeM},
\cite{Lev}. 
\enddemo
\msk
\example{6.6. Gromoll filtration} The Gromoll
filtration of $x\in \pi_0\wdt{\DIFF}(\DD^{i-1})$ is the largest number
$j=j(x)$  such that $x$ lifts from 
$$\pi_0\wdt{\DIFF}(\DD^{i-1})\cong\pi_{j-1}\wdt{\DIFF}(\DD^{i-j})$$
to $\pi_{j-1}\DIFF(\DD^{i-j})$. This is 
the original definition of \cite{Grom}; see also \cite{Hit}. 

To obtain upper bounds on  $j(x)$ in special cases, we use 6.5, with
$k=i+1$ and $n=i-j$ and $M=\DD^n$.
Therefore: if $4$ divides $i+1$
and $x$ has Gromoll filtration $\ge j$, and is the image of 
$\bar x\in  L_{i+1}(\ZZ)$, then the
signature of $\bar x$ is divisible by $2a_{j+1}$ (by $4a_{j+1}$ in the case
where $4$ divides $j+1$). 
\endexample
\msk
\example{6.7. Riemann--Roch for smooth fiber bundles}
Let $p:E\to B$ be a fiber bundle with fibers $\cong M$ 
and structure group $\DIFF(M,\partial M)$.
Let $R$ be a ring, and let $V$ be a bundle (with discrete 
structure group) of f\.g\. left proj\. $R$--modules on $E$. 
This determines $[V]:E\to\Omega\8\bK(R)$. Let 
$V_i$ be the bundle on $B$ with fiber $H_i(p^{-1}(b);V)$ over $b$. 
We assume that the fibers of $V_i$ are projective; then each $V_i$
determines $[V_i]:B\to K(R)$. Then the following Riemann--Roch formula holds:
$$\text{tr}^*[V]=\sum(-1)^i[V_i]\qquad\qquad\in[B,\Omega\8\bK(R)]
\tag6.7.1$$
where $\text{tr}:\Sigma\8B_+\to \Sigma\8E_+$
is the Becker--Gottlieb--Dold transfer \cite{BeGo}, \cite{Do}, \cite{DoP}, a stable map 
determined by $p$ . Both sides of \thetag{6.7.1} have meaning
when $p:E\to B$  is a fibration whose fibers are homotopy 
equivalent to compact CW--spaces. However, \thetag{6.7.1} does 
not hold in this generality. It can fail for a fiber bundle
with compact (and even closed) {\it topological} manifolds as fibers. 

Formula \thetag{6.7.1} is a distant corollary of 5.2.4.  Namely, 
both \thetag{6.7.1} and 5.2.4 are ways of saying that 
certain generalized Euler characteristics of 
a smooth  compact $M$ lift canonically to 
$\Omega\8\Sigma\8(M_+)$. For the proof of \thetag{6.7.1}, see \cite{DWW}. 
Earlier, Bismut and Lott \cite{BiLo} had proved by 
analytic methods that \thetag{6.7.1} holds in the case $R=\CC$ after
certain  characteristic classes are applied to both sides of the 
equation. 
\endexample
\msk
\example{6.8. Obstructions to finding block bundle structures} Suppose
that 
$p:E\to B$  is a fibration with connected base whose fibers are
oriented Poincar\'e duality spaces  of formal dimension $2k$. Let
$E\td\to E$ be a normal covering  with translation group $\pi$. 
With these data we can associate a map 
$$B\la \Omega\8\bhm_{\pi}(k)\tag6.8.1$$
where $\bhm_{\pi}(k)$ is the (topological, connective) $K$--theory 
of f\.g\. projective $\RR\pi$--modules with nondegenerate 
$(-1)^k$--hermitian form \cite{Wa1}, \cite{Wa4}. The map \thetag{6.8.1} 
stably classifies the hermitian bundle on $B$ with fiber 
$H^k(E\td_x;\RR)$ over $x\in B$. There is a {\it hyperbolic map} \cite{Wa4}
from $\bo_{\pi}$, the (topological, connective) $K$--theory 
of f\.g\. projective $\RR\pi$--modules, to $\bhm_{\pi}(k)$. 
Let $\bhm_{\pi}(k)/\bo_{\pi}$ be its mapping cone. The map  
$$B\la\Omega\8(\bhm_{\pi}(k)/\bo_{\pi})\tag6.8.2$$
obtained by composing \thetag{6.8.1} with $\Omega\8$ of 
$\bhm_{\pi}(k)\to \bhm_{\pi}(k)/\bo_{\pi}$ is the {\it family
multisignature} of $p$. Now suppose that 
$p$ admits a block bundle structure; see \S5.1. Then \thetag{6.8.2} factors {\it
rationally} as 
$$B\la \Omega\8(\bhm_e(k)/\bo_e)@>\quad\otimes\RR\pi\quad>>
\Omega\8(\bhm_{\pi}(k)/\bo_{\pi})$$ 
where $e$ is the trivial group.  
The case $B=*$ appears in
\cite{Wa1, \S13B}. The general statement can be proved by 
expressing the rationalized family multisignature in terms of the
family  visible symmetric signature of \S5.1.
\endexample
\msk
\example{6.9. Relative calculations of $h$--structure Spaces} Using \S1.6 and 
4.2.1 one obtains, in the case where $M$ is smooth, a diagram 
$$\sS_d(M^n)@>\quad\nabla\quad>>\sS_d(\tau)@>{\quad \lambda\quad}>>
\Omega^{\infty+n}\bL\bA\bul^h(M,n)
\tag{6.9.1}$$
which is a homotopy fiber sequence in the concordance 
stable range;  more precisely, the composite map in \thetag{6.9.1} is trivial and 
the resulting map from $\sS_d(M)$ to the homotopy fiber of $\lambda$ is 
approximately $(n-1)/3$--connected. This formulation has some relative
variants which are attractive   because relative $LA$--theory is often
easier  to describe than absolute $LA$--theory, and the estimates 
available for the relative concordance stable range 
are often better than those for the absolute concordance 
stable range. 

{\it Illustration:}  Let $M$ be smooth, compact, connected, with
connected boundary, and suppose the inclusion 
$\partial M\to M$ induces an isomorphism of fundamental groups. Let
$M_0=M\minus\partial M$.  Define $\sS_d(M_0)$ as the Space
of pairs $(N,f)$ where 
$f:N\to M_0$ is a proper homotopy equivalence of smooth manifolds 
without boundary.  Let $\sS_d(\tau_0)$ be the corresponding Space of 
$n$--dimensional vector bundles on $M_0$ equipped with a stable fiber 
homotopy equivalence to $\tau_0:=\tau|M_0$. Using some controlled 
$L$--theory and algebraic $K$--theory of spaces,  one
obtains the following variation  on \thetag{6.9.1}: a diagram 
$$
\sS_d(M_0)@>\quad\nabla\quad>>\sS_d(\tau_0)@>\qquad>>
\Omega^{\infty+n}(\bL\bA\bul^h(M,n)/\bL\bA\bul^h(\partial M,n))
\tag6.9.2$$
which is a homotopy fiber sequence in the 
$\le (n/3-c_1)$ range, where $c_1$ is a constant independent of $n$.
Our hypothesis on fundamental groups implies that 
$\bL\bA\bul^h(M,n)/\bL\bA\bul^h(\partial M,n)
\simeq(\SS^n_!\wedge(\bA(M)/\bA(\partial M)))\ho$
and \cite{BHM},  \cite{Go3} imply 
$$\bA(M)/\bA(\partial M)\,\simeq\, \bT\bC(M)/\bT\bC(\partial M)$$ 
where $\bT\bC$ is the 
{\it topological cyclic homology} viewed as a spectrum--valued functor. See also 
\cite{Ma1},  \cite{Ma2}, \cite{HeMa}. 
So \thetag{6.9.2} simplifies to 
$$\sS_d(M_0)@>\quad\nabla\quad>>\sS_d(\tau_0)@>\qquad>>
\Omega^{\infty+n}((\SS^n_!\wedge(\bT\bC(M)/\bT\bC(\partial
M)))\ho)\,.$$ This is still a homotopy fiber sequence in the $\le
(n/3-c_1)$ range. Unpublished  results of  T.Goodwillie and G.Meng
indicate that it is a fibration sequence in the 
$\le n-c_2$ range for some constant $c_2$,  provided
$\partial M\to M$ is $2$--connected.
\endexample
\msk

\enddocument